\documentclass[a4paper, 10pt]{article}  \usepackage{amsmath}
\usepackage{amssymb}
\usepackage{amsthm}       \newtheorem{thm}{Theorem}[section]
\newtheorem{prop}[thm]{Proposition}
\newtheorem{lem}[thm]{Lemma}
\newtheorem{cor}[thm]{Corollary}  \theoremstyle{definition}
\newtheorem{df}[thm]{Definition}   \theoremstyle{definition}
\newtheorem{rem}{Remark}                \theoremstyle{plain}
 \theoremstyle{definition}
\newtheorem{ex}[thm]{Example}   \def\CC{\Bbb{C}}
\def\RR{\Bbb{R}}  
\def\CCI{\overline{\CC}}        \def\NN{\Bbb{N}} 
\def\B1{{\rm\kern.32em\vrule    width.12em       height1.4ex
depth-.05ex\kern-.28em 1}}

\begin{document}
\title{Semi-hyperbolic fibered rational maps and rational semigroups
\footnote{published in Ergodic Theory and Dynamical Systems (2006), 
{\bf 26}, 893-922.}}
 \author{Hiroki Sumi\
\\   Department of Mathematics,\\ Graduate School of Science,\\ 
Osaka University, \\ 1-1, \ Machikaneyama,\ Toyonaka,\ Osaka,\ 560-0043,\ 
Japan\\ E-mail:    sumi@math.sci.osaka-u.ac.jp\\ 
http://www.math.sci.osaka-u.ac.jp/$\sim $sumi/welcomeou-e.html
}   
\maketitle
\begin{abstract}
This paper is based on the author's 
previous work \cite{S4}. 
We consider fiber-preserving 
complex dynamics on fiber bundles whose fibers are 
Riemann spheres and whose base spaces are 
compact metric spaces. 
In this context, without any assumption on 
(semi-)hyperbolicity, we show that the fiberwise Julia sets 
are uniformly perfect. 
From this result, we show that, 
for any semigroup $G$ generated by 
a compact family of rational maps on 
the Riemann sphere 
$\CCI $ of degree two or greater,\ 
the Julia set of any 
subsemigroup of $G$ is uniformly perfect.

 We define the semi-hyperbolicity of dynamics on fiber bundles and
 show that, if the dynamics on a fiber bundle 
is semi-hyperbolic,\ then the fiberwise 
Julia sets are porous, and 
the dynamics is weakly rigid. 
Moreover, we show that if the dynamics is semi-hyperbolic and 
the fiberwise maps are polynomials,\ then under some conditions, 
the fiberwise basins of infinity are John domains.
We also show that the Julia set of 
a rational semigroup (a semigroup generated 
by rational maps on $\CCI $) that is 
semi-hyperbolic, except at perhaps 
finitely many points in the Julia set, and 
which satisfies  
the open set condition, is either porous or 
equal to the closure of the open set.
Furthermore, we derive an upper estimate of the Hausdorff dimension 
of the Julia set.  
\end{abstract}
\section{Introduction}
\label{Introduction}
 To investigate random one-dimensional 
complex dynamics,\ the dynamics of 
semigroups generated by rational maps on 
the Riemann sphere $\CCI $ and  
fiber-preserving holomorphic dynamics on 
fiber bundles that appear in complex dynamics  
in several dimensions,\ we consider the  
dynamics of fibered rational maps;\ that is,\ 
we consider fiber-preserving complex dynamical systems on fiber bundles 
whose fibers are Riemann spheres and whose base 
spaces are general compact metric spaces.
 The notion of the dynamics of fibered rational maps,\ which 
 is a generalized notion of the ``dynamics of the fibered 
 polynomial maps" presented
 by O. Sester(\cite{Se1}, \cite{Se2}, \cite{Se3}),\ 
  was
 introduced by M. Jonsson in \cite{J2}. Researches 
on the dynamics of semigroups generated by rational maps on 
the Riemann sphere (\cite{HM1},\ \cite{HM2},\ 
\cite{HM3}\, \cite{GR},\ \cite{Bo},\ 
 \cite{St1},\ \cite{St2},\ \cite{St3},\ 
 \cite{SSS},\ \cite{S1},\ \cite{S2},\ \cite{S3},\   
\cite{S4},\ \cite{S5},\ \cite{S6},\ \cite{S7}),\ the random iterations of rational
 functions(\cite{FS},\ \cite{BBR},\ \cite{Br})
  and polynomial skew products 
 on $\CC ^{2}$ (\cite{H1},\ \cite{H2},\ \cite{J1}) 
 are directly related to this subject. For 
 polynomial skew products (the dynamics of fibered polynomials)
 whose base spaces are 
 general compact metric 
 spaces,\ 
 see O.Sester's work \cite{Se1}, \cite{Se2} and \cite{Se3}.
  In \cite{Se3}, he investigated the quadratic case in detail. 
  In particular,\ 
  he developed a combinatorial theory for quadratic fibered 
  polynomials and  constructed an abstract space of 
  combinatorics. Moreover, he demonstrated the realizability and 
  rigidity of abstract combinatorics.  

  For the ergodic theory of random diffeomorphisms,\  
 see \cite{K}. 

 For the Hausdorff dimension of the Julia set of 
 a (semi-) hyperbolic rational semigroup, 
 see \cite{S4} and \cite{S7}.

 For dynamics of semi-hyperbolic transcendental semigroups, 
 see \cite{KS1}, \cite{KS2}.

  In this paper, by applying some of the results in \cite{S4} 
  obtained by the author,\ we show that 
 semi-hyperbolicity along the fibers of a fibered rational map 
  implies that the fiberwise Julia sets are $k$-porous, 
  where the constant $k$ does not depend on any 
  points of the base space (hence the upper box dimensions 
  of the fiberwise Julia sets are uniformly less than a constant 
  that is strictly less than $2$)
  ({\bf Theorem~\ref{porositythm}});\ 
  that the dynamics have a kind of weak rigidity 
({\bf Theorem~\ref{rigidity}}); and that 
if the fiberwise maps are polynomials with some conditions, then 
the fiberwise basins of infinity are $c$-John domains, 
where $c$ is a constant not depending on any points 
in the base spaces({\bf Theorem~\ref{semihypJohn}}). 
This generalizes a result from \cite{CJY} and a result from 
\cite{Se3}. These results are stated 
in Section~\ref{Results on f}, and the 
proofs are given in Section~\ref{Proofs on f}.
For research on the 
semi-hyperbolicity of 
the usual dynamics of rational functions,\ 
see \cite{CJY} and \cite{Ma}. 

 In Section~\ref{Results on r},\ we 
 present several results on the dynamics of 
 rational semigroups; i.e., on the semigroups 
 generated by rational maps on the Riemann sphere.
We show that the Julia set of a finitely generated 
 rational semigroup that is semi-hyperbolic, 
 except at most finitely many points in the 
 Julia set, and that satisfies the ``(backward) open set condition"
 (which is a generalized notion in research 
 on self-similar sets constructed by some similarity transformations) 
 with 
 an open set,\ 
 is either porous (and hence the upper box dimension 
 is strictly less than $2$) or equal to the closure 
 of the open set ({\bf Theorem~\ref{rsporosity}}). 
 Furthermore, under the same assumptions, we show that 
 the Hausdorff dimension of the Julia set 
 is less than the infimum of numbers $\delta $, which 
 allows the construction of a ``$\delta $-subconformal measure"\ 
 that is less than the critical exponent of the semigroup
 ({\bf Theorem~\ref{uhfindim}}).    
 The proofs of these results are given in 
 Section~\ref{Proofs on r}.

 A proof of the above results requires the 
 potential theoretic stories (Section~\ref{Potential Theory}),\ 
 distortion lemmas for holomorphic proper maps,\  
 various results on non-constant limit functions, 
 and the continuity of the fiberwise Julia sets
 with respect to points in the base spaces
 (this is an important property and is not easy to show), 
 from \cite{S4}. Note that the upper semicontinuity 
 of fiberwise Julia sets with respect to points in the 
 base spaces does not hold in general.
 
 Furthermore, in Section~\ref{Results on f}, we show that 
 the fiberwise Julia sets of a fibered rational map 
 are $C_{1}$-uniformly perfect 
 and that the Hausdorff dimensions are greater than 
 a positive constant $C_{2}$,
 where $C_{1}$ and $C_{2}$ are   
 constants that do not depend on any points in the base 
 space,\ if the fiberwise maps are of degree 
 two or greater,\ even if semi-hyperbolicity is not 
assumed ({\bf Theorem~\ref{unifperf1}}).
 The proof is given in Section~\ref{Proofs on f}. 
  This result on uniform perfectness is used to 
  show the result on Johnness. Furthermore, 
  from the uniform perfectness of fiberwise 
  Julia sets,\ we show that, for any semigroup $G$ generated by 
a compact family of rational maps on 
$\CCI $ of degree two or greater,\ 
there exists a positive constant 
$C$ such that the Julia set of any 
subsemigroup of $G$ is $C$-uniformly perfect
({\bf Theorem~\ref{unifperrs}}).
 The proof,\ which is given in Section~\ref{Proofs on r},\  
 is based on combining 
 the density of repelling fixed points in the Julia set 
 (an application of Ahlfors's five-island theorem) 
 with the potential theory in Section~\ref{Potential Theory}.

\subsection{Results on fibered rational maps}
\label{Results on f} 
In this section, we state the results on fibered 
rational maps. The proofs are given in 
Section~\ref{Proofs on f}. 
First, we provide some notation and definitions regarding 
the dynamics of fibered rational maps.
\begin{df}( \cite{J2}) 
\label{cbundledf}
A triplet $(\pi , Y, X)$ is called a ``$\CCI $-bundle" 
if 
\begin{enumerate}
\item $Y$ and $X$ are compact metric spaces, 
\item $\pi :Y\rightarrow X$ is a continuous and surjective map, 
\item \label{localtriv}
       and there exists an open covering $\{ U_{i}\} $  
       of $X$ such that, for each $i $, there exists a 
       homeomorphism  $\Phi _{i}:U_{i}\times \CCI 
       \rightarrow \pi ^{-1}(U_{i})$ such  that
       $ \Phi _{i}(\{ x\} \times \CCI )=\pi ^{-1}\{ x\} $ 
       and $\Phi _{j}^{-1}\circ \Phi _{i}:
       \{ x\} \times \CCI \rightarrow  
       \{ x\} \times \CCI $ is a M\"{o}bius map 
       for each $x\in U_{i}\cap U_{j},\ $ 
       under the identification 
       $\{ x\} \times \CCI  \cong \CCI .$     
\end{enumerate}
\begin{rem}
\label{frsdfrem}
 By condition \ref{localtriv},\ each fiber 
$Y_{x}:=\pi ^{-1}\{ x\} $ has a complex structure.
Furthermore,  
given $x_{0}\in X$, one may find a continuous 
family $i_{x}:\CCI \rightarrow Y_{x}$ of homeomorphisms, 
for $x$ close to $x_{0}.$  Such a family $\{ i_{x}\} $ 
is called a ``local parameterization". 
Since $X$ is compact,\ we may assume throughout
this paper that there exists a compact subset 
$M_{0}$ of the set of M\"{o}bius transformations 
of $\CCI $ such that $i_{x}\circ j_{x}^{-1}\in M_{0}$ 
for any two local parameterizations $\{ i_{x}\} $ and 
$\{ j_{x}\} $ .  

 Moreover, throughout this paper, we assume the 
 following condition:
 \begin{itemize}
 \item  there exists a smooth $(1,1)$-form 
$\omega _{x}>0$ inducing the distance on $Y_{x}$ from $Y$, and 
$x\mapsto \omega _{x}$ is continuous.  That is, 
if $\{ i_{x}\} $ is a local parameterization,\ 
then the pull back $i_{x}^{\ast }\omega _{x}$ is 
a positive smooth form on $\CCI $ that depends continuously 
on $x$ with respect to the $C^{\infty }$ topology.

\end{itemize}
\end{rem}
\end{df}
\noindent {\bf Notation:} If 
$(\pi , Y=X\times \CCI ,\ X)$ is a trivial $\CCI $-bundle,\ 
then we denote by $\pi _{\CCI }:Y\rightarrow \CCI $ 
the second projection.
\begin{df}
Let $(\pi ,Y,X)$ be a $\CCI $-bundle.
Let $f:Y\rightarrow Y$ and $g:X\rightarrow X$ be 
continuous maps. Let $f$ be called a fibered rational map 
 over $g$ (or a rational map fibered over $g$), if 
\begin{enumerate}
\item $\pi \circ f=g\circ \pi $
\item $f|_{Y_{x}}:Y_{x}\rightarrow Y_{g(x)}$ is a rational map, 
for any $x\in X.$ That is, 
$(i_{g_{x}})^{-1}\circ f\circ i_{x} $ 
is a rational map from $\CCI $ to itself, for any 
 local parameterization $i_{x}$ at $x\in X$ and 
 $i_{g(x)}$ at $g(x).$  
\end{enumerate}
\end{df}

\noindent {\bf Notation: } If $f:Y\rightarrow Y$ is a fibered rational map 
over $g:X\rightarrow X, $ then we set   
$f_{x}^{n}=f^{n}|_{Y_{x}}:Y_{x}\rightarrow Y_{g^{n}(x)}$, for any $x\in X$ and $n\in \NN .$ Moreover, we set 
$f_{x}= f|_{Y_{x}}:Y_{x}\rightarrow Y_{g(x)}.$ 
Furthermore, we set $d_{n}(x)=\deg (f_{x}^{n})$ and 
$d(x)=d_{1}(x)$, for any $x\in X $ and $n\in \NN .$ 
Note that the definition of fibered rational map implies that 
the function $x\mapsto d(x)$ is continuous on $X.$
\begin{df}
Let $(\pi , Y,X)$ be a $\CCI $-bundle. 
Let $f:Y\rightarrow Y$ be a fibered rational map 
over $g:X\rightarrow X. $ Then, for any $x\in X$, we denote
by $F_{x}(f)$(or simply $F_{x}$) the set of points $y\in Y_{x}$ 
that has a neighborhood $U$ of $y$ in $Y_{x}$ such that 
$\{ f_{x}^{n}\} _{n\in \NN }$ is a normal family in $U;$
that is, $y\in F_{x}$ if and only if 
the family $Q_{x}^{n}=i_{x_{n}}^{-1}\circ f_{x}^{n}\circ 
i_{x} $ of rational maps on $\CCI $ 
($x_{n}:=g^{n}(x)$ ) is normal near 
$i_{x}^{-1}(y)$. Note that, by 
Remark~\ref{frsdfrem},\ 
 this does not depend on the choices of local parameterizations 
at $x$ and $x_{n}.$ Equivalently, $F_{x}$ is the 
open subset of $Y_{x}$, where the family $\{ f_{x}^{n} \} $ 
of mappings from $Y_{x}$ into $Y$ is locally equicontinuous.
We set $J_{x}(f)$(or simply $J_{x}$)$=Y_{x}\setminus F_{x}.$  

 Furthermore, we set 
$ \tilde{J}(f)=\overline{\bigcup _{x\in X}J_{x}} $, 
where the closure is taken in the space $Y,\ $  
$\tilde{F}(f)=Y\setminus \tilde{J}(f),$ 
and $\hat{J}_{x}(f)$(or simply $\hat{J}_{x}$)$=
\tilde{J}(f)\cap Y_{x} $, for each $x\in X.$  
\end{df}
\begin{rem}
There exists a fibered rational map $f:Y\rightarrow Y$ 
such that $\bigcup _{x\in X}J_{x}$ is 
NOT compact. 

\end{rem}

\begin{df}
\label{bundleex1}
(\cite{S4})
Let $h_{1},\ldots ,h_{m}$ be non-constant 
rational maps. Let 
$\Sigma _{m}=\{ 1,\ldots ,m\} ^{\NN }$ 
be the space of one-sided infinite sequences of 
$m$ symbols, and 
$\sigma :\Sigma _{m}\rightarrow \Sigma _{m}$ 
be the shift map; that is,\ 
$\sigma $ is defined by 
$\sigma ((w_{1},w_{2},\ldots ))=(w_{2},w_{3},\ldots ).$ 
Let $X$ be a compact subset of $\Sigma _{m}$ 
such that $\sigma (X)\subset X.$ 
Let $Y=X\times \CCI $ and $\pi :Y\rightarrow X$ 
be the natural projection. 
Then $(\pi ,Y,X)$ is a $\CCI $-bundle.
Let $f:Y\rightarrow Y$ be a map defined by:
$f((w,y))=(\sigma (w),\ h_{w_{1}}(y))$, 
where $w=(w_{1},w_{2},\ldots ).$      
 Then, 
$f:Y\rightarrow Y$ is a fibered rational map 
over $\sigma :X\rightarrow X.$ 

In the above, if $X=\Sigma _{m}$, 
then  
$f:Y\rightarrow Y$ is called the 
{\bf fibered rational map 
associated with the generator system 
$\{ h_{1},\ldots h_{m}\} $ } (of 
the rational semigroup $G=\langle h_{1},\ldots ,h_{m}\rangle  $). 
(See Section~\ref{Results on r} for 
the definition of rational semigroups.) 
\end{df}
Now we present a result on uniform perfectness.
\begin{df}
Let $C$ be a positive number. Let 
$K$ be a closed subset of $\CCI .$ We say that 
$K$ is {\bf $C$-uniformly perfect} if, 
$\sharp K\geq 2$ and 
for any doubly connected domain $A$ in $\CCI $ such 
that $A$ separates $K$; i.e., such that 
$A\subset \CCI \setminus K$ and 
both of the two connected components of $\CCI \setminus A$ 
have non-empty intersections with $K, $ 
 mod $A$ (the modulus 
of $A$; for the definition,\ see \cite{LV}) is less than $C.$  
\end{df}
\begin{thm}({\bf Uniform perfectness})
\label{unifperf1}
Let $(\pi , Y ,  X)$ be a  
$\CCI $-bundle. Let $f:Y\rightarrow Y$ be a 
fibered rational map over $g:X\rightarrow X$ 
such that $d(x)\geq 2$,  for any $x\in X.$ 
(In this theorem, we do NOT assume semi-hyperbolicity.)
Then, we have the following.
\begin{enumerate}
\item 
\label{unifperf1-1}
There exists a
 positive constant $C$ 
 such that, for any $x\in X,\ $ 
 $J_{x}$ and $\hat{J}_{x}$ are $C$-uniformly 
 perfect. Furthermore, there exists a positive constant
  $C_{1}$ such that for each $x\in X$,\ we have diam $J_{x}\geq C_{1}$, 
   with respect to the distance 
  in $Y_{x}$ induced by $\omega _{x}.$ Moreover,\ 
  there exists a positive constant $C_{2}$ such that, 
  for each $x\in X$, the Hausdorff dimension $\dim _{H}(J_{x})$ of 
  $J_{x}$, with respect to the distance 
  on $Y_{x}$ induced by $\omega _{x}$, satisfies the 
  condition that $\dim _{H}(J_{x})\geq C_{2}.$ (Note that 
  $C_{1}$ and $C_{2}$ do not depend on $x\in X.$ ) 
\item 
\label{unifperf1-2}
Suppose further that $f(\tilde{F}(f))\subset \tilde{F}(f)$(
for example,\ assume that $g:X\rightarrow X$ 
is an open map).
If a point $z\in Y$ satisfies 
$f_{\pi (z)}^{n}(z)=z$ and $(f_{\pi (z)}^{n})'(z)=0$ 
for some $n\in \NN $ and $z\in \hat{J}_{\pi (z)},\ $ 
then $z$ belongs to the interior of $\hat{J}_{\pi (z)}$ 
with respect to the topology of $Y_{\pi (z)}.$ 

 \end{enumerate}
\end{thm}
Theorem~\ref{unifperf1} is now used to show 
a result on Johnness (Theorem~\ref{semihypJohn}) and 
Theorem~\ref{unifperrs}. 
\begin{ex}
Let $z_{0}\in \CCI $ be a point. 
Let $h_{1}$ and $h_{2}$ be two rational maps on $\CCI $ 
of degree two or greater. 
Let $f:\Sigma _{2}\times \CCI \rightarrow \Sigma _{2}\times \CCI $ 
be the fibered rational map associated with the 
generator system $\{ h_{1},h_{2}\} .$ 
Suppose that $z_{0}$ is a superattracting fixed point of $h_{1}$ 
and is a repelling fixed point of $h_{2}.$ Then, it can  easily be seen 
that $z_{0}\in \hat{J} _{x}$, where $x=(1,1,\ldots )\in \Sigma _{2}.$ 
Since the shift map $\sigma :\Sigma _{2} \rightarrow \Sigma _{2}$ 
is an open map,\ by Theorem~\ref{unifperf1} it follows that $z_{0}$ belongs to 
the interior of $\hat{J} _{x}.$ 
\end{ex}
\begin{rem}
Theorem~\ref{unifperf1} generalizes Theorem 3.1 in \cite{Br}. 
\end{rem}
\begin{rem}
\label{uprem}
Uniform perfectness implies many useful properties
(\cite{BP},\cite{Po},\cite{Su}). This terminology was introduced in \cite{Po}.
For a survey on uniform perfectness, see \cite{Su}. 
We now consider the following:
\begin{enumerate}
\item In \cite{BP}, it was shown that 
a closed subset $K$ of $\CCI $ with $\sharp K\geq 3$ is $C$-uniformly perfect
if and only if there exists a constant $\delta $ such that, 
for any component $U$ of $\CCI \setminus K,\ $ 
\begin{equation}
\label{updfeq1}
\lambda _{U}(z)>\delta /\mbox{dist}(z,\partial U),\ 
\end{equation} 
where $\lambda _{U}(z)$ denotes the density of the 
hyperbolic metric of $U$ at $z$ and dist$(z,\partial U)$ 
denotes the Euclidian distance of the point $z$ 
from the set $\partial U.$ (If $K$ is bounded and $U$ 
is the unbounded component of $\CCI \setminus K,\ $ then
inequality (\ref{updfeq1}) will hold only for all $z\in U$ sufficiently 
close to $K.$) In the above discussion, $\delta $ depends only on $C$, and 
$C$ depends only on $\delta .$ This fact will be used to 
show Theorem~\ref{semihypJohn}.
Detailed inequalities regarding the relationships among 
$\delta ,\ C$ and other invariants are presented in \cite{Su}. 
  
\item If a closed subset $K$ of $\CCI $ is 
$C$-uniformly perfect, then 
the Hausdorff dimension $\dim _{H}(K)$ of 
$K$ with respect to the spherical metric 
satisfies $\dim _{H}(K)\geq C'>0,\ $ 
where $C'$ is a positive constant that depends 
only on $C$ (Theorem 7.2 in \cite{Su}).
      
\end{enumerate}
\end{rem}

\noindent {\bf Notation :} 
\begin{itemize}
\item 
Let $Z_{1}$ and $Z_{2}$ be two topological spaces and 
$g: Z_{1}\rightarrow Z_{2}$ be a map. For any 
subset $A$ of $Z_{2},\ $ We denote by $c(A,g)$ the set 
of all connected components of $g^{-1}(A).$  
\item For any $y\in \CCI $ and $\delta >0, $ we set 
$B(y,\delta )=\{ y'\in \CCI \mid d(y,y')<\delta \} ,$ 
where $d$ is the spherical distance. Similarly,\ 
for any $y\in \CC $ and $\delta >0 $, 
we set $D(y,\delta )=\{ y'\in \CC \mid |y-y'|<\delta \} .$  
\item Let $(\pi ,Y,X)$ be a $\CCI $-bundle. 
For any $y\in Y$ and $r>0$, we set 
$ \tilde{B}(y,r)=\{ y'\in Y_{\pi (y)}\mid d_{\pi (y)}(y',y)<r\} ,\ $
where, for each $x\in X$, 
we denote by $d_{x}$ the distance on $Y_{x}$ induced by the form $\omega _{x}.$

\end{itemize}
We now define the semi-hyperbolicity of 
fibered rational 
maps. 
\begin{df}({\bf semi-hyperbolicity of fibered rational maps})
Let\\ $(\pi ,Y,X)$ be a $\CCI $-bundle. Let 
$f:Y\rightarrow Y$ be a fibered rational map over
 $g:X\rightarrow X.$ Let $N\in \NN .$ 
We denote by $SH_{N}(f)$ the set of points $z\in Y$ 
 such that there exists a positive number $\delta, $
 a neighborhood $U$ of $\pi (z)$, and a local parameterization 
 $\{ i_{x}\} $ in $U$ such that, 
 for any $x\in U,$ any $n\in \NN ,$ 
  any $x_{n}\in g^{-n}(x)$ and any $V\in c(i_{x}
  (B(i_{\pi (z)}^{-1}(z),\ \delta )),\ f_{x_{n}}^{n}),\ $ 
  we have that 
  $$ \deg (f_{x_{n}}^{n}:V\rightarrow 
  i_{x}(B(i_{\pi (z)}^{-1}(z),\ \delta )))
  \leq N.$$  We set 
   $UH(f)=Y\setminus \bigcup _{N\in \NN }SH_{N}(f).$ 
   A point $z\in SH_{N}(f)$ is called a 
   {\bf semi-hyperbolic 
   point of degree }$N.$
   Furthermore, we say that $f$ is {\bf semi-hyperbolic (along fibers) } 
   if $\tilde{J}(f)\subset \bigcup _{N\in \NN }SH_{N}(f).$
   This is equivalent to $\tilde{J}(f)\subset SH_{N}(f)$, 
   for some $N\in \NN .$ 
\end{df}
\begin{df}
Let $f:Y\rightarrow Y$ be a rational 
map fibered over $g:X\rightarrow X.$ 
Then, we set 
$C(f):=\bigcup _{x\in X}\{ \mbox{critical points of }f_{x}\}\  
(\subset Y). $ This is called the {\bf fiber critical set } of 
fibered rational map $f.$ Furthermore,\  
we set 
$ P(f)=\overline{\bigcup _{n\in \NN }\bigcup _{x\in X}
f_{x}^{n}(\mbox{ critical points of }f_{x})}$ where 
the closure is taken in $Y.$ 
This is called the {\bf fiber postcritical set} of the fibered 
rational map $f.$ 
Furthermore, we say that $f:Y\rightarrow Y$ is 
{\bf hyperbolic 
along fibers} if 
 $P(f)\subset F(f). $ (Then $f$ is semi-hyperbolic 
along fibers with the constant $N=1.$) 

\end{df}
We now present a result on the John condition. 

\noindent {\bf Notation:} Let 
$y\in \CC $ and $b\in \CCI $ be two distinct points. 
Furthermore, let $E$ be a curve in $\CCI $ joining $y$ to $b$ 
such that $E\setminus \{ b\} \subset \CC .$ 
Then, for any $c\geq 1$, we set 
$$ \mbox{car }(E,c,y,b)=\bigcup _{z\in E\setminus 
\{ y,b\} }D(z,\ \frac{|y-z|}{c}).$$ 
This is called the $c$-carrot with core $E$ and 
vertex $y$ joining $y$ to $b.$   
\begin{df}
Let $V$ be a subdomain of $\CCI $ such that 
$\infty \in \CCI \setminus \partial V.$  Let 
$c\geq 1$ be a number. Then, we say that 
$V$ is a {\bf $c$-John domain} if there exists a 
point $y_{0}\in \overline{V}$ such that,  
for any $y\in V\setminus \{ y_{0}\} $, 
there exists a curve $E$ joining $y_{0}$ to 
$y$ such that $E\setminus \{ y_{0}\} \subset \CC $ 
and 
$ \mbox{car }(E,c,y,y_{0})\subset V.$ 
In what precedes, the point $y_{0}$ is called the 
center of the John domain $V.$  
\end{df}
\begin{df}
Let $(\pi ,\ Y=X\times \CCI ,\ X)$ be a trivial $\CCI $-bundle.
A fibered rational map $f:Y\rightarrow Y$ over 
$g:X\rightarrow X$ is called a 
{\bf fibered polynomial map over $g:X\rightarrow X$} 
if $f$ is of the form:
$ f(x,y)=(g(x),q_{x}(y))$ for 
$(x,y)\in X\times \CCI $,\ where 
$q_{x}$ is a polynomial for each $x\in X.$ 
For such an $f$,\ we set 
$A_{x}(f):=\{ y\in Y_{x}\mid 
\pi _{\CCI }(f_{x}^{n}(y))\rightarrow \infty \} .$  
\end{df}
\begin{thm}({\bf Johnness})
\label{semihypJohn}
Let $(\pi, Y=X\times \CCI, \ X)$ be a trivial
$\CCI $-bundle. Let $f:Y\rightarrow Y$ be a
semi-hyperbolic fibered polynomial map over $g:X\rightarrow X$
such that $d(x)\geq 2$, for any $x\in X.$
Assume that either 
\begin{enumerate}
\item 
 $J_{x}(f)$ is connected for each $x\in X$,\ 
or 
\item  \label{semihypJohnh2}
the map $x\mapsto f_{x}$ from
$X$ to the space of polynomials is locally constant (here, we
identify $f_{x}$ with a polynomial on $\CCI $) and \\ 
$\inf \{ d(a,b)\mid a\in \pi _{\CCI }(J_{x}),\ 
b\in \pi _{\CCI }(A_{x}(f)\cap C(f)),\ x\in X\} >0.$ 
\end{enumerate}
Then, there exists a
positive constant $c$ such that, for any $x\in X$,
the basin of
infinity $A_{x}(f)$
is a $c$-John domain.
\end{thm}
\begin{rem}
\begin{itemize}
\item 
For a fibered polynomial map $f$ such that 
$d(x)\geq 2$ for each $x\in X$,\ 
we have that 
$J_{x}(f)$ is connected for each $x\in X$ if and only if 
$\pi _{\CCI }(P(f))\setminus \{ \infty \} $ is bounded in 
$\CC .$ 
\item 
In Theorem~\ref{semihypJohn},
if $X$ is a set consisting of one point, \ then 
$f$ is semi-hyperbolic if and only if 
the basin of infinity is a John domain(\cite{CJY}). 
In \cite{Se3}, O. Sester made the same assertion 
as that presented in Theorem~\ref{semihypJohn} for 
``non-recurrent quadratic fibered polynomials" 
with connected fiberwise Julia sets.  
\end{itemize}
\end{rem}

\begin{rem}
Johnness likewise implies many useful properties
(\cite{NV},\ \cite{Jone}).For example,\ 
if $V$ is a John domain,\ then 
the following facts hold.
\begin{itemize}
\item 
If $\infty \in V,\ $ then 
the center of $V$ is $\infty .$ 
\item 
Let $a\in \partial V $ 
and $b\in V.$ Then there exists a curve 
$E$ joining $a$ to $b$ and a constant $c$ 
such that car $(E,c,a,b)\subset V.$ In 
particular,\ 
$a$ is accessible from $b.$ 
\item 
$V$ is finitely connected at any point in $\partial V$:
\ that is,\ if $y\in \partial V,\ $ then there exists 
an arbitrarily small open neighborhood $U$ of $y$ in $\CCI $ 
such that $U\cap V$ has only finitely many connected components.    
\item If $V$ is simply connected, 
then it follows that $\partial V$ is locally connected. 
\item 
$\partial V$ is 
holomorphically removable:\ that is,\ if 
$\varphi :\CCI \rightarrow \CCI $ is a homeomorphism 
and is holomorphic on $\CCI \setminus \partial V,\ $ then 
$\varphi $ is holomorphic on $\CCI .$ From this fact,\ 
 we deduce that the two-dimensional Lebesgue measure 
of $\partial V$ is equal to zero.

\end{itemize}
\end{rem} 
In order to present results on porosity and other properties,\  
let us first establish some technical conditions. 
\begin{df}[Condition(C1)]
Let $(\pi ,Y,X)$ be a $\CCI $-bundle.
 Let $f:Y \rightarrow Y$  
be a fibered rational map over $g:X\rightarrow X.$ 
 We say that $f $ 
satisfies condition (C1) if there exists a family  
 $\{ D_{x}\} _{x\in X}$ of topological disks with 
 $D_{x}\subset Y_{x},\ x\in X$ such that the 
following conditions are satisfied: 
\begin{enumerate}
\item for each $x\in X$, there exists 
a point $z_{x}\in Y_{x}$ and a positive 
number $r_{x}$ such that 
$D_{x}=\tilde{B}(z_{x},r_{x}),$ 
\item $\overline{\bigcup _{x\in X}\bigcup _{n\geq 0}f_{x}^{n}( D_{x})}\subset 
\tilde{F}(f), $  
\item for any  $x\in X ,\ $ it holds that 
diam$(f_{x}^{n}(D_{x}))\rightarrow 0,\ $ as 
$n\rightarrow \infty ,$ and 
\item $\inf _{x\in X} r_{x}>0.$
\end{enumerate} 
\end{df} 
\begin{df}[Condition(C2)]
Let $(\pi ,Y,X)$ be a $\CCI $-bundle.
 Let $f:Y \rightarrow Y $ 
be a fibered rational map over $g:X\rightarrow X.$ 
 We say that $f $ 
satisfies condition (C2) if, for each $x_{0}\in X $, 
there exists an open neighborhood $O$ of $x_{0}$
 and a family $\{ D_{x}\} _{x\in O}$ of topological disks 
 with $D_{x}\subset Y_{x}, x\in O$  
such that the following conditions are satisfied:   
\begin{enumerate}
\item for each $x\in O$, there exists 
a point $z_{x}\in Y_{x}$ and a positive 
number $r_{x}$ such that 
$D_{x}=\tilde{B}(z_{x},r_{x}),$ 

\item $\overline{\bigcup _{x\in O}\bigcup _{n\geq 0}f_{x}^{n}( D_{x})}\subset 
\tilde{F}(f), $  
\item for any  $x\in O ,\ $ 
it holds that diam$(f_{x}^{n}(D_{x}))\rightarrow 0,\ $ as 
$n\rightarrow \infty ,$ and 
\item $x\mapsto D_{x}$ is continuous in $O.$ 
\end{enumerate} 
\end{df}
\begin{rem}
\begin{enumerate}
\item Let $(\pi ,\ Y=X\times \CCI ,\ X)$ be 
a trivial $\CCI $-bundle. 
Let $f:Y\rightarrow Y$ be a 
fibered polynomial map over $g:X\rightarrow X$ such 
that $d(x)\geq 2$ for each $x\in X.$   
Then, setting $D_{x}=D$, where $D$ is a small neighborhood 
of infinity, for each $x\in X,\ $ 
the fibered rational map $f$ satisfies condition (C2) 
with the family of disks $(D_{x})_{x\in X}.$  
\item 
Let $\{ h_{1},\ldots h_{m}\} $ be non-constant 
rational functions on $\CCI $ with $\deg (h_{1})\geq 2.$ 
Let $f:Y\rightarrow Y$ be the fibered rational map 
associated with the generator system 
$\{ h_{1},\ldots ,h_{m}\} $ of rational semigroup 
$G=\langle h_{1},\ldots ,h_{m}\rangle $
(See Section~\ref{Results on r}). 
Suppose that $f$ is semi-hyperbolic along fibers 
and that $\pi _{\CCI }(\tilde{J}(f))=J(G)$
(See Section~\ref{Results on r} and \cite{S5}) 
 is 
not equal to the Riemann sphere. 
Moreover, suppose that each element of Aut $(\CCI ) \cap G$
( if this is not empty ) is loxodromic.
Then, 
it follows that 
$f$ satisfies condition (C2). 
Actually,\ there exists an attracting fixed  
point $a$ of some element of $G$ in $F(G).$ Since 
$G$ is semi-hyperbolic,\  
if $D_{x}=D(a,\epsilon)$, for each $x\in \Sigma _{m}$ 
where $\epsilon $ is a positive number,\ 
then $f$ satisfies condition (C2) with the 
family of disks $(D_{x})_{x\in \Sigma _{m}}.$    
For more details,\ see Theorem 1.35 and Remark 5 in \cite{S4}.
\end{enumerate}
\end{rem}
We now present a result on porosity.
\begin{df}
\label{porositydf}
Let $(Y,d)$ be a metric space. 
Let $k$ be a constant such that 
$0<k<1.$ Let $J$ be a 
subset of $Y.$ We say that 
$J$ is {\bf $k$-porous} if, 
for any $x\in J$ and 
any positive number $r$, there exists 
a ball in $\{ y\in Y\mid d(y,x)<r\} \setminus J$ with 
radius at least $kr.$  
\end{df}
\begin{rem}
\label{porosityrem}
If $Y$ is the Euclidean space $\RR ^{n}$ 
and $d$ 
is the Euclidean metric,\ the  
upper box dimension of any $k$-porous 
bounded set $J$ in $\RR ^{n}$ 
is less than 
$n-c(k,n),\ $ where 
$c(k,n)$ is a positive constant 
that depends only on $k$ and $n$
 (Proposition 4.2 and the proof in \cite{PR}).
\end{rem}
\begin{thm}({\bf Porosity})
\label{porositythm}
Let $(\pi ,Y,X)$ be a $\CCI $-bundle. 
Let $f:Y\rightarrow Y$ be a fibered rational map 
over $g:X\rightarrow X.$ 
Suppose that $f$ satisfies condition (C1)
and that $f$ is semi-hyperbolic. 
Then there exists a constant $0<k<1$
such that, for each $x\in X$,\ $J_{x}$ is $k$-porous 
in $Y_{x}.$ In particular,\ 
there exists a constant $0\leq c<2$ such that, 
for each $x\in X,\ $ 
$ \dim _{H}(J_{x})\leq \overline{\dim }_{B}(J_{x})\leq c,\ $
where $\dim _{H}$ denotes the Hausdorff dimension 
and $\overline{\dim }_{B}$ denotes the upper box dimension with respect 
to the metric on $Y_{x}$ induced by $\omega _{x}$
($\omega _{x}$ is of the form in Remark~\ref{frsdfrem}).  

\end{thm}
We now present a sufficient condition 
for $\tilde{J}=\bigcup _{x\in X}J_{x}$ and 
the complement of the
backward images of non-semi-hyperbolic points 
in fiberwise Julia sets to be of measure zero.
\begin{df}({\bf Good points for fibered rational maps})
\label{goodf}
Let $(\pi , Y,X)$ be a $\CCI $-bundle. 
Let $f:Y\rightarrow Y$ be a fibered rational map 
over $g:X\rightarrow X.$ 
We set 
$ \tilde{J}_{good}(f)=
\{ z\in \tilde{J}(f)\mid \limsup _{n\rightarrow \infty }
d(f^{n}(z),UH(f))>0\} .$ 
\end{df}
\begin{thm}({\bf Measure zero})
\label{bundlemeas}
Let $(\pi ,Y,X)$ be a $\CCI $-bundle. 
Let $f:Y\rightarrow Y$ be a fibered rational map 
over $g:X\rightarrow X.$ Suppose 
that all of the following conditions hold:
\begin{enumerate}
\item $f$ satisfies condition (C1);
\item for each $x\in X,\ $ 
the boundary of $\hat{J}_{x}(f)\cap UH(f)$ 
in $Y_{x}$ does not separate points in 
$Y_{x}$
(i.e. $Y_{x}\setminus \partial (\hat{J}_{x}(f)\cap UH(f))$ 
is connected);  
\item $\tilde{J}(f)
\setminus 
\bigcup _{n\in N}f^{-n}(UH(f))\subset \tilde{J}_{good}(f); $
 and 
\item \label{noshrink}
for each $z\in \tilde{J}(f)\cap UH(f)$ 
and each open neighborhood $V$ of $z$ 
in $Y_{\pi (z)}$, 
the diameter of $f^{n}_{\pi (z)}(V)$ does not tend to zero 
as $n\rightarrow \infty .$ 
\end{enumerate} 
Then, $\tilde{J}(f)=\bigcup _{x\in X}J_{x}$ and, 
for each 
$x\in X,\ $ the two-dimensional Lebesgue 
measure of $J_{x}\setminus \bigcup _{n\in N}f^{-n}(UH(f))$ 
is equal to zero.
\end{thm} 
The statement $\tilde{J}(f)=\bigcup _{x\in X}J_{x}$
in Theorem~\ref{bundlemeas} is used to show Theorem~\ref{rsporosity}.

  We now present a result on rigidity.
\begin{thm}({\bf A rigidity})
\label{rigidity}
Let $(\pi , Y,X)$ and $(\tilde{\pi } , \tilde{Y}, \tilde{X})$ 
be two $\CCI $-bundles. Let $f:Y\rightarrow Y$ be 
a fibered rational map over $g:X\rightarrow X$ 
and $\tilde{f}:\tilde{Y}\rightarrow \tilde{Y}$ 
be a fibered rational map over $\tilde{g}:\tilde{X}
\rightarrow \tilde{X}.$  Let $u:Y\rightarrow \tilde{Y}$ 
be a homeomorphism that is a bundle conjugacy between 
$f$ and $\tilde{f}$ ; i.e., $u$ satisfies the condition that 
$\tilde{\pi }\circ u=v\circ \pi $ for some homeomorphism 
$v:X\rightarrow X$ and $\tilde{f}\circ u=u\circ f.$ 
Suppose that $f$ is semi-hyperbolic along fibers 
and satisfies condition (C1).  For each $w\in X,\ $
let $u_{w}:Y_{w}\rightarrow 
\tilde{Y}_{v(w)}$ be the restriction of $u.$ 
Let $x\in X$ be a 
point and assume that 
$u_{x}$ is $K$-quasiconformal on $F_{x}.$ Then,\  
for each 
 $a\in \overline{\cup _{n\in \Bbb{Z}} g^{n}(\{ x\} )}$, 
it follows that  $u_{a}:Y_{a}\rightarrow \tilde{Y}_{v(a)}$ 
 is $K$-quasiconformal on the whole $Y_{a}$,\ with 
 the same dilatation constant $K.$     

\end{thm}

\subsection{Results on rational semigroups}
\label{Results on r}
In this section, we present several results on 
rational semigroups. The proofs are given 
in Section~\ref{Proofs on r}.
 Before stating results, we will first establish some notation and definitions regarding the  dynamics of rational semigroups.

 For a Riemann surface  $S$, let End\(  (S) \) denote the set
of all holomorphic endomorphisms of $S$. 
It is a semigroup whose semigroup operation is 
the functional composition. A {\bf rational semigroup} is a 
subsemigroup of End($\CCI $) without  any constant  elements.    
We say that 
a rational semigroup  $G$  is a  {\bf   polynomial  semigroup} 
if  each element of $G$ is a polynomial. Research 
on the dynamics of  rational semigroups was initiated by 
A. Hinkkanen and G. J. Martin (\cite{HM1}),\ 
who were interested in the role that
the dynamics of polynomial semigroups plays in research 
on various one-complex-dimensional 
moduli spaces for discrete groups,\  
 and by
F. Ren's group(\cite{GR}).

 In this paper,\ a rational semigroup $G$ may not be 
 finitely generated,\ in general.   
\begin{df}

Let $G$ be a rational semigroup. We set
\[ F(G) = \{ z\in \CCI \mid G \mbox{ is normal in a neighborhood of  $z$} \} ,
\ J(G)  = \CCI \setminus  F(G) .\] \(  F(G)\) is  called the
{\bf Fatou set}  of  $G$, and \( J(G)\)  is  called the {\bf 
Julia set} of $G$. The backward orbit $G^{-1}(z)$ of $z$ and 
the 
{\bf set of exceptional points } $E(G)$ are defined\
 by:
 $ G^{-1}(z) =\cup _{g\in G}  g^{-1}\{ z\}  $ and 
$ E(G)= \{z\in \CCI \mid \sharp G^{-1}(z) \leq2 \} .$
Furthermore,\ we set 
$G(z):= \cup _{g\in G}\{ g(z)\} .$ 
For any subset $A$ of $\CCI ,\ $ set 
$G^{-1}(A)=\cup _{g\in G}g^{-1}(A).$  
We denote by $\langle h_{1},h_{2},\ldots \rangle $ the 
rational semigroup generated by the family $\{ h_{i}\} .$ 
For a rational map $g$, we denote by $J(g)$ the Julia set 
of the dynamics of $g.$
\end{df}
We now present a result on uniform perfectness. In the following 
theorem,\ we do NOT assume semi-hyperbolicity.
\begin{thm}({\bf Uniform perfectness})
\label{unifperrs}
Let $\Lambda $ be a compact set 
in the space $\{ h:\CCI \rightarrow \CCI \mid 
h:\mbox{holomorphic,\ }\deg (h)\geq 2 \} $ 
endowed with 
topology induced by uniform convergence on $\CCI .$ 
Let $G$ be a rational semigroup
 generated by the set $\Lambda .$ Then 
 there exists a positive constant $C$ such that 
  each subsemigroup $H$ of $G$ satisfies the condition that 
  $J(H)$ is $C$-uniformly perfect.   
Furthermore, for any subsemigroup $H$ of $G$, 
if a point $z_{0}\in J(H)$ satisfies the condition that
there exists an element $h\in H$ such that 
$h(z_{0})=z_{0}$ and $h'(z_{0})=0$, then 
it follows that $z_{0}\in $int $J(H).$ 
\end{thm} 
\begin{ex}
Let $G=\langle h_{1},h_{2}\rangle $ 
where $h_{1}(z)= 2z+z^{2}$ and 
$h_{2}(z)=\frac{z^{3}}{z-a}$,\  
$a\in \CC $ with $a\neq 0.$ 
Then, $0\in J(G)$ and $0$ is a 
superattracting fixed point of 
$h_{2}.$ Hence, 
$0\in $int $J(G)$,\ by 
Theorem~\ref{unifperrs}. 
Since $\infty $ is a common 
attracting fixed point of $h_{1}$ and 
$h_{2}$, we have 
$\infty \in F(G).$ 
Furthermore, 
let $H$ be a subsemigroup of $G$ such that 
$0\in J(H)$ and 
$h_{2}\in H.$ Then, by Theorem~\ref{unifperrs} again, 
we have $0\in $int $J(H).$ Moreover, 
we have $\infty \in F(H).$ 

 In particular, let $H_{0}$ be a 
 subsemigroup of $G$ that is generated by 
 $G\setminus \langle h_{1}\rangle .$ 
 Then we have all of the following: 
 \begin{enumerate}
 \item $0\in $ int $J(H_{0}).$
 \item $\infty \in F(H_{0}).$ 
 \item For any finitely generated 
 subsemigroup $H_{1}$ of 
 $H_{0}$, we have $0\in F(H_{1}).$
 \end{enumerate}
 For, since $0\in $int $J(G)$, 
 $J(G)=\overline{\cup _{g\in G}J(g)}$ 
 (Corollary 3.1 in \cite{HM1}), 
 and $J(h_{1})$ is nowhere dense, 
 we obtain that there exists a 
 suquence $(g_{n})$ in $H_{0}$ such that 
 $d(0,J(g_{n}))\rightarrow 0$ as $n\rightarrow \infty .$ 
 Hence, $0\in  J(H_{0}).$ Since 
 $0$ is a superattracting fixed point of 
 $h_{2}$ and $h_{2}\in H_{0}$, 
 by Theorem~\ref{unifperrs} we obtain 
 $0\in $ int $J(H_{0}).$ 
 Since $H_{0}\subset G$ and $\infty \in F(G)$, 
 we obtain $\infty \in F(H_{0}).$ 
 For any finitely generated subsemigroup 
 $H_{1}$ of $H_{0}$, since $0$ is a common 
 attracting fixed point of any element of $H_{0}$, 
 it follows that $0\in F(H_{1}).$  
 \end{ex} 

\begin{rem}
In \cite{HM2}, by A. Hinkkanen and G. Martin, 
it was shown that a finitely generated 
rational semigroup $G$ such that 
each $g\in G$ is of degree two or greater  
satisfies the condition that $J(G)$ is uniformly perfect. 
In \cite{St3}, by R. Stankewitz, it was shown that, 
if $\Lambda $ (this is allowed to have an 
element of degree one) is a 
family of rational maps on $\CCI $ such that 
the Lipschitz constant of each element of 
$\Lambda $ with respect to the spherical metric 
on $\CCI $ is uniformly bounded,\ then 
 the Julia set of semigroup $G$ generated by 
$\Lambda $ is uniformly perfect. However, there 
has been no research done on the uniform perfectness 
of the Julia sets of subsemigroups of such a semigroup.
In \cite{HM2} and \cite{St3}, the proofs were based on 
the density of repelling fixed points in the 
Julia sets (which was shown by an 
application of Ahlfors's five-island theorem),\ 
whereas, in this paper, the proof of 
Theorem~\ref{unifperrs} is based on 
the combination of the density of repelling fixed points 
in the Julia set with Proposition~\ref{potentialprop}
(potential theory).   
\end{rem}

We now define the semi-hyperbolicity of rational semigroups.
 \begin{df}({\bf Semi-hyperbolicity of rational semigroups})
\label{shndf}
Let  $G$  be  a rational semigroup  and  $N$ a  positive
integer. We denote by $SH_{N}(G)$ 
the set of points $z\in \CCI  $ such that 
there exists a positive number $\delta $ 
such that, 
for any $g\in G$ and any $V\in c(B(z,\delta ),\ g),\ $ 
it holds that 
$ \deg (g:V\rightarrow B(z,\delta ))\leq N.$  
Furthermore, we set $UH(G)=\CCI \setminus (\cup _{N\in \NN }SH_{N}(G)).$
A point $z\in SH_{N}(G)$ is called a 
{\bf semi-hyperbolic point of degree }$N$ (for $G$).
We say that $G$  is {\bf 
 semi-hyperbolic  }   if
$J(G)\subset
 \bigcup _{N\in \NN }SH_{N}(G).$ 
 This is equivalent to $J(G)\subset SH_{N}(G)$, for some $N\in \NN .$

 Furthermore,\ for a rational semigroup $G$,\ we set \\ 
 $P(G):=\overline{\bigcup _{g\in G}\{ 
 \mbox{critical values of } g\} }\ (\subset \CCI ).$ 
 This is called the {\bf postcritical set } of $G.$ 
 We say that $G$ is {\bf hyperbolic } if 
 $P(G)\subset F(G).$ 
 Note that if $G$ is hyperbolic, then $G$ is semi-hyperbolic.
Note also that if $G=\langle h_{1},\ldots ,h_{m}\rangle$,\ 
then $P(G)=\overline{\bigcup _{g\in G}
g\left( \bigcup _{j=1}^{m}\{ \mbox{ critical values of } h_{j}\}\right) } .$  
\end{df}
\begin{df}
\label{oscdf}
Let $G =  \langle  h_{1},h_{2},\ldots , h_{m} \rangle$  be  a
finitely generated rational semigroup. 
Let $U$ be a non-empty open set in $\CCI .$ 
 We say  that  $G$
satisfies the  {\bf open set  condition} 
with $U$ with respect to the
generator systems  $\{ h_{1},h_{2},\ldots ,h_{m}\} $,   if 
  for  each   $j=1,\ldots   ,  m,\
h_{j}^{-1}(U)\subset U$ and 
$\{ h_{j}^{-1}(U)\} _{j=1,\ldots ,m} $ are mutually disjoint.
\end{df}
We now present a result on porosity.
\begin{thm}({\bf Porosity})
\label{rsporosity}
Let $G=\langle h_{1},\ldots ,h_{m}\rangle $ 
be a rational semigroup with an element of 
degree two or greater. 
Suppose that all of the following conditions hold:
\begin{enumerate}

\item $G$ satisfies the open set condition 
with an open set $U$ with respect to the 
generator system $\{ f_{1},\ldots ,f_{m}\} ,\ $

\item $\sharp (UH(G)\cap J(G))<\infty $ and 

\item $UH(G)\cap J(G)\subset U.$
\end{enumerate} 

Then it follows that $J(G)=\overline{U}$, 
or that $J(G)$ is porous (and so the upper box dimension of 
$J(G)$ is strictly less than $2$). 
Moreover,\ the fibered rational map 
$f:\Sigma _{m}\times \CCI \rightarrow \Sigma _{m}\times \CCI $ 
associated with the generator system $\{ h_{1},\ldots ,h_{m}\} $ 
satisfies 
$ \tilde{J}(f)=\bigcup _{x\in \Sigma _{m}}J_{x}.$ 
\end{thm}
We now present results regarding the Hausdorff dimension of 
Julia sets.
\begin{df}
Let $G$ be a rational semigroup and $\delta $ a non-negative number. 
We say that a Borel probability measure $\mu $ on $\CCI $ is 
{\bf $\delta $-subconformal} (for $G$) if, for each $g\in G$ and for each Borel-measurable set $A$, it holds that 
$ \mu (g(A))\leq \int _{A} \| g'(z)\| ^{\delta }d\mu ,$
where we denote by $\| \cdot \| $ the norm of the derivative 
with respect to the spherical metric. 
We set 
$s(G)=\inf \{\delta \mid \exists \mu : 
\delta \mbox{-subconformal measure}\}  .$ 

 Let $G$ be a countable rational semigroup.  
For each $x\in \CCI $ and each real number $s$, we set
$ S(s, x)=\sum _{g\in G}\sum _{g(y)=x}\|g'(y)\| ^{-s } $,
counting multiplicities and 
$ S(x)=\inf \{ s\mid S(s, x)<\infty \} .$ 
If there is no $s$ such that $S(s, x)<\infty ,$\ then we set $S(x)=\infty .$ 
Also,\ we set 
$ s_{0}(G)=\inf \{  S(x)\mid x\in \CCI \}.$ 
 
\end{df}
An estimate on $s_{0}(G)$ exists when $G$ satisfies the open set condition.
\begin{prop}
\label{critexp}
Let $G=\langle h_{1},\ldots ,h_{m} \rangle $ 
be a rational semigroup. When $m=1,\ $ 
assume that $h_{1}$ is neither the identity nor 
an elliptic M\"{o}bius transformation. 
Suppose that $G$ satisfies the open set condition 
with an open set $U$ with respect to the 
generator system $\{ h_{1},\ldots ,h_{m}\}.$ 
Suppose also that $J(G)\neq \overline{U}.$ 
Then there exists an open set $V$ included in 
$U\cap F(G)$ such that, for almost every $x\in V$ with respect 
to the two-dimensional Lebesgue measure,\  
$S(2,x)<\infty .$ 
In particular,\ $s_{0}(G)\leq 2.$  
\end{prop}
\begin{thm}({\bf Hausdorff dimension})
\label{uhfindim}
Let $G=\langle h_{1},\ldots ,h_{m} \rangle $ 
be a rational semigroup. When $m=1$,\ 
we assume $J(G)\neq \CCI .$ 
Under the same assumptions as those of 
Theorem~\ref{rsporosity},\ it follows 
that 
$\dim _{H}(J(G))\leq s(G)\leq s_{0}(G),\ $
where $\dim _{H}$ denotes the Hausdorff dimension 
with respect to the spherical metric in $\CCI .$  
\end{thm}
\begin{ex}
Let $h_{1}(z)=z^{2}+2,\ h_{2}(z)=z^{2}-2$ 
and $U=\{ |z|<2.\} .$ Then, 
$h_{1}^{-1}(U)\cup h_{2}^{-1}(U)\subset U$ 
and $h_{1}^{-1}(U)\cap h_{2}^{-1}(U)=\emptyset .$
Let $h_{3}$ be a 
polynomial that is conjugate to $h_{4}^{n}$ 
by an affine map $\alpha ,\ $ where 
$h_{4}(z)=z^{2}+\frac{1}{4}$ and $n\in \NN$ 
is a sufficiently large number. For an appropriate $\alpha $,\ 
it holds that $J(h_{3})\subset U\setminus 
(h_{1}^{-1}(\overline{U})\cup h_{2}^{-1}(\overline{U})).$ 
For sufficiently large $n$,\ 
$h_{3}^{-1}(\overline{U})\subset U\setminus 
(h_{1}^{-1}(\overline{U})\cup h_{2}^{-1}(\overline{U})).$
Then $G=\langle h_{1},h_{2},h_{3}\rangle $ 
satisfies the conditions specified as assumptions in 
Theorem~\ref{rsporosity}. In this case, 
$UH(G)\cap J(G)$ is one point that is 
a parabolic fixed point $a$ of 
$h_{3}$ (For, $UH(G)\cap J(G)$ is 
included in $P(G)\cap J(G) 
=\{ -2,2,a\} $, and, by Lemma~\ref{rsporolem0.5},
it follows that $\{ -2,2\} \subset \CCI \setminus (UH(G)\cap J(G))$) 
and $UH(G)\cap J(G)\subset U.$
 By Theorem~\ref{rsporosity},\ 
it can be seen that $J(G)$ is porous and that, in particular,\ 
the upper box dimension is strictly less than two. 
Furthermore, by Theorem~\ref{uhfindim} we obtain that 
$\dim _{H}(J(G))\leq s(G)\leq s_{0}(G).$
\end{ex}
\subsection{Remark}
 We now make some remarks about fibered rational maps 
and rational semigroups.
\begin{rem}
\begin{enumerate}
\item Dabija \cite{D} showed that (almost) every holomorphic 
 selfmap $f$ of a ruled surface $Y$ is a fibered rational map over a 
 holomorphic map $g:X\rightarrow X, $ where $X$ is a Riemann 
 surface.
\item \label{bundleex1-3}
Let $p(x)\in \CC [x]$ be a polynomial 
of degree two or greater, and 
$q(x,y)\in \CC [x,y]$ be a polynomial 
of the form: $q(x,y)=y^{n}+a_{1}(x)y^{n-1}+\cdots .$ 
Let $f:\CC ^{2}\rightarrow \CC ^{2}$ 
be a map defined by 
$ f((x,y))=(p(x),q(x,y)).$ 
This is called a polynomial skew product in $\CC ^{2}.$ 
  Dynamics of maps of this form were investigated by S.-M. Heinemann 
  in \cite{H1} and  \cite{H2}, and by 
  M. Jonsson in \cite{J1}. 
Let $X$ be a compact subset of $\CC $ such that 
$p(X)\subset X.$ (e.g., the Julia set of $p.$) 
Let $(\pi, Y=X\times \CCI ,X)$ 
be a trivial $\CCI $-bundle.
Then the map $\tilde{f}:Y\rightarrow Y$ 
defined by $\tilde{f}((x,y))=(p(x), q(x,y))$ 
is a fibered rational map over 
$p:X\rightarrow X.$ 
\item 
\label{bundlesemiex1.5}
In Corollary 6.7 of \cite{Se3}, O. Sester showed that 
any ``non-recurrent quadratic fibered polynomials" with 
connected fiberwise Julia sets are semi-hyperbolic.  
\item \label{bundlesemiex2}
Let $\{ h_{1},\ldots ,h_{m}\} $ be 
 non-constant rational functions on $\CCI .$ 
 Let $f:Y\rightarrow Y$ be the 
 skew product map associated with 
 the generator system $\{ h_{1},\ldots ,h_{m}\} .$ 
 It can be shown by means of simple arguments that 
 $f:Y\rightarrow Y$ is semi-hyperbolic along 
 fibers if and only if 
 $G=\langle h_{1},\ldots ,h_{m}\rangle $ is semi-hyperbolic.
\item In \cite{S4},\ with $G$ a finitely generated rational semigroup,\ 
 a sufficient condition for the semi-hyperbolicity 
 of a point $z\in J(G)$ was given,\ which gives a 
 generalization of R. Ma\~{n}\'{e}'s work(\cite{Ma}). 
 Furthermore, in \cite{S4}, a necessary and sufficient condition 
 for the semi-hyperbolicity of a finitely generated 
 rational semigroup was given. From this, it can be seen that 
 $G=\langle z^{2}+2,z^{2}-2 \rangle $ is semi-hyperbolic,\ 
 and is not hyperbolic.    
\end{enumerate}

\end{rem}
\section{Tools and Proofs}
We now present the proofs of the main results, 
meanwhile providing further notation and tools. 
\subsection{Fundamental properties of fibered rational maps}
By means of definitions,\ the following lemma can be easily shown.
\begin{lem}
\label{frlem}
Let $(\pi ,Y,X)$ be a $\CCI $-bundle. 
Let $f:Y\rightarrow Y$ be a fibered rational 
map over $g:X\rightarrow X.$ Then, 
\begin{enumerate}
\item 
\label{frlem1}
For each $x\in X,\ $ $f_{x}^{-1}(F_{g(x)})=F_{x},\ 
f_{x}^{-1}(J_{g(x)})=J_{x}.$ Furthermore,\ 
$f(\tilde{J}(f))\subset \tilde{J}(f).$
\item 
\label{frlem2}
If $g:X\rightarrow X$ is an open 
map,\ then 
$f^{-1}(\tilde{J}(f))=\tilde{J}(f)$ and 
$f(\tilde{F}(f))\subset \tilde{F}(f).$
\item 
\label{frlem3}
If $g:X\rightarrow X$ is a surjective 
and open map,\ then 
$f^{-1}(\tilde{J}(f))=\tilde{J}(f)=f(\tilde{J}(f))$ 
and $f^{-1}(\tilde{F}(f))=\tilde{F}(f)=f(\tilde{F}(f)).$
\end{enumerate}
\end{lem}
\begin{proof}
This proof is the same as that for Lemma 2.4 in
\cite{S4}.
\end{proof}
\subsection{Fundamental properties of 
rational semigroups}
For a rational semigroup $G$, for each $f\in G$, it holds that 
$f(F(G))\subset F(G)$ and $f^{-1}(J(G))\subset J(G).$ Note that  
equality $f^{-1}(J(G))= J(G)$ does not hold in general.
If \( \sharp J(G)\geq 3 \) ,\ then \( J(G) \) is a 
perfect set,\  
$ \sharp E(G) \leq 2,\ $  
\( J(G) \) is the smallest closed backward invariant set containing at least
three points, and $J(G)$ is the closure of the union of all 
repelling fixed points of the elements of $G$,\  
which implies that 
$J(G)=\overline{\bigcup _{g\in G}
J(g)}.$ 
If a point \( z\) is not in \( E(G),\ \) then, for every \( x\in J(G),\ \)
$x\in   \overline{G^{-1}(z)} .\) In particular, 
if \ \( z\in  J(G)\setminus E(G) \) \ then
$ \overline{G^{-1}(z)}=J(G).$
For more precise statements,\ see 
Lemma 2.3 in \cite{S5},\ for which the proof  
is based on \cite{HM1} and \cite{GR}. 
Furthermore,\ if $G$ is generated by a precompact
 subset $\Lambda $ of End$(\CCI )$, then 
$ J(G)=\overline{\bigcup _{f\in \Lambda }f^{-1}(J(G))}=
\bigcup _{h\in \overline{\Lambda }}h^{-1}(J(G)).$ 
In particular, if $\Lambda $ is compact, then 
$J(G)= \bigcup _{f\in \Lambda }f^{-1}(J(G))$(\cite{S4}).
 We call this property of the Julia set 
 {\bf backward self-similarity}. 
\begin{rem}
In the context of backward self-similarity,\ the existing research 
on the Julia sets of rational semigroups 
may be considered as a kind of generalization 
of the research on self-similar sets constructed 
by some similarity transformations from $\CC  $ to itself,\ 
which can be regarded as the Julia sets of some 
rational semigroups. It can be easily seen that 
the Sierpi\'{n}ski gasket is the Julia set 
of a rational semigroup 
$G=\langle h_{1},h_{2},h_{3}\rangle $, where 
$h_{i}(z)=2(z-p_{i})+p_{i},i=1,2,3$, with 
$p_{1}p_{2}p_{3}$ a regular triangle. 
\end{rem}
\subsection{Potential theory and measure theory}
\label{Potential Theory}
For the proof of results on uniform perfectness, 
Johnness, and so on,  
let us borrow some notation from \cite{J2} and \cite{S4},\ 
concerning potential theoretic aspects. 
By the arguments in \cite{J2} 
and \cite{S4}, for a fibered 
rational map $f:Y\rightarrow Y$ 
over $g:X\rightarrow X$ 
with 
$d(x)\geq 2$ for each $x\in X$, 
one can show a result corresponding 
to Proposition 2.5 in \cite{S4},
using the arguments in $\S 3$ in \cite{J2} 
and from pp. 580-581 in \cite{S4} . 
In this paper, the following statements,\  
and especially the lower semicontinuity of 
$x\mapsto J_{x}(f)$, are necessary.
(Proposition~\ref{potentialprop}.\ref{potentialprop11}).
\begin{prop}
\label{potentialprop}
Let $(\pi ,Y,X)$ be a $\CCI $-bundle.
Let $f:Y \rightarrow Y $ be a rational 
map fibered over $g:X\rightarrow X. $ 
 Assume that  
$d(x)\geq 2$, for each $x\in X.$ Then, 
for each $x\in X$, there exists a
 Borel probability measure 
 $\mu _{x}$ on $Y$ satisfying all of the following.
\begin{enumerate}
\item \label{potentialprop8}
$x\mapsto \mu _{x}$ is continuous with respect to the weak topology 
of probability measures in $Y .$ 
\item \label{potentialprop9}
supp$(\mu _{x}) =J_{x}$, for each $x\in X.$ 
\item \label{potentialprop11}
$x\mapsto J_{x}$ is lower semicontinuous with respect to 
the Hausdorff metric in the space of 
the non-empty compact subsets of $Y .$  
That is,\ if $x,x_{n}\in X, x_{n}\rightarrow x$ as 
$n\rightarrow \infty  $ and $y\in J_{x},\ $ then 
there exists a sequence $(y_{n})$ of points in $Y$ with 
 $y_{n}\in J_{x_{n}}$, for each $n\in \NN $, such that 
 $y_{n}\rightarrow y$ as $n\rightarrow \infty .$     
\end{enumerate} 
Furthermore , $J_{x}(f)$ is a non-empty perfect set, for 
 each $x\in X.$
\end{prop}
 \begin{proof}
 Since $d(x)\geq 2$ for each $x\in X$, and 
 $x\mapsto d(x)$ is continuous,\ 
 we can demonstrate these statements in the same way as 
  in \S 3 in \cite{J2}, using the argument
 from pp. 580-581 in \cite{S4}. 
 The statement \ref{potentialprop11} follows 
 easily from \ref{potentialprop8} and \ref{potentialprop9}. 
  
 \end{proof}
\subsection{Distortion lemma and 
limit functions}
For an exposition of the results deduced from 
semi-hyperbolicity, the following 
statements are necessary.
For the research on the semi-hyperbolicity of 
the usual dynamics of rational functions,\ 
see \cite{CJY} and \cite{Ma}.

\noindent {\bf Notation:} 
\begin{enumerate}
\item 
Let $X$ be a compact set in $\CCI $ and 
$z$ be a point in $\CCI \setminus X.$ Then 
we set
$ \mbox{Dist}(X,z)=\max\limits _{y\in X}d(y,z)/
\min\limits _{y\in X}d(y,z).$
\item For two positive numbers $A$ and $B,\ $ 
$A\asymp B$ means $K^{-1}\leq A/B\leq K$, for some 
constant $K$ independent of $A$ and $B.$   
\end{enumerate} 

\begin{lem}[\cite{CJY}]
\label{ylem1} ({\bf Distortion lemma for proper maps})
For any positive integer $N$ and real number $r$ with $0< r
< 1,\ $ there exists a constant $C=C(N,r)$ such that, if $ f:
D(0,1)\rightarrow D(0,1)$ is  a proper holomorphic map  with
$\deg (f) =N $, \  then
$ H(f(z_{0}), C)\subset f(H(z_{0}, r))\subset H(f(z_{0}),r)$ 
for any $z_{0}\in  D(0,1)$,\ where 
$H(x,s)$ denotes the hyperbolic disk in $D(0,1)$ 
around 
$x$ with radius $s.$  Here, the constant $C=C(N,r)$
depends only on $N$ and $r.$
\end{lem}
The following is a generalized distortion lemma for proper maps.
\begin{lem}[\cite{S4},\cite{S6}]
\label{backsmall}
Let $V$ be a domain in $\CCI ,\ $ $K$, a continuum in $\CCI $ with 
diam$_{S}K=a.$  Assume $V\subset \CCI \setminus K.$ Let $f:V\rightarrow D(0,1)$ be a proper holomorphic map of degree $N.$ Then 
there exists a constant $r(N,a)$ depending only on $N$ and $a$ such that, for each $r$ with $0<r\leq r(N,a),\ $ there exists a constant $C=C(N,r)$ depending only on $N$ and $r$ such that, for each connected component $U$ of $f^{-1}(D(0,r)),\ $ it follows that 
$ \mbox{diam}_{S}\ U\leq C,$
where we denote by diam$_{S}$ the spherical diameter.    
Furthermore, it follows that $C(N,r)\rightarrow 0$ as $r\rightarrow 0.$ 
\end{lem}
As the author mentioned in
Remark 6 in \cite{S4},\ 
the results Lemma 2.13, Theorem 2.14, Lemma 2.15 and Theorem 2.17 
in \cite{S4} are generalized to the version of 
the dynamics of fibered rational maps on 
$\CCI $-bundles,\ which we need to 
show results in this paper.
The following lemma is a slightly 
modified version of 
Lemma 2.15 in \cite{S4}.
\begin{lem}[\cite{S4}]
\label{keythmlem}
Let $(\pi, Y,X)$ be a $\CCI $-bundle. 
Let $f:Y\rightarrow Y $ be a 
fibered rational map over $g:X\rightarrow X.$  
Assume that $f$ satisfies condition (C1). Assume 
$z_{0}\in SH_{N}(f)$, for some $N\in \NN .$ 
Then  there
 exists a positive number $\delta _{0}$ such that, for each 
$\delta $ with $0<\delta <\delta _{0}$,  
there exists a neighborhood $U$ of $x_{0}:=\pi (z_{0})$ 
in $X $ such that, for each $n\in \NN , $ each $x\in U $ and 
each $x_{n}\in p^{-n}(x),\ $ it follows that each element of 
$c(i_{x}i_{x_{0}}^{-1}\tilde{B}(z_{0}, \delta ),\ f_{x_{n}}^{n})$ 
is simply connected.  
\end{lem}
The following theorem relates to what happens 
if there exists a non-constant limit function on 
a component of a fiber-Fatou set.
This is the key to state other results.
The proof is the same as that for 
Lemma 2.13 in \cite{S4}. 
\begin{thm}[\cite{S4}]
\label{keyth}({\bf Key theorem I})
Let $(\pi ,Y,X)$ be a $\CCI $-bundle.
Let $f:Y\rightarrow Y$ be a fibered rational map 
over $g:X\rightarrow X.$ 
Assume that $f$ satisfies condition (C1). 
Let $z\in Y$ be a point with 
$z\in F_{\pi (z)}.$ Let 
$(i_{x})$ be a local parameterization.
Let $U$ be a connected open neighborhood of 
$i_{\pi (z)}^{-1}(z)$ in $\CCI .$ 
Suppose that there exists a 
strictly increasing 
sequence 
$(n_{j})$ of $\NN $ such that 
$R_{j}:=i_{\pi f^{n_{j}}(z)}^{-1}\circ f_{\pi (z)}^{n_{j}}\circ 
i_{\pi (z)}$ uniformly converges to a non-constant map $\phi $  
on $U$ as $j\rightarrow \infty .$ Furthermore, suppose that 
$f_{\pi (z)}^{n_{j}}(z)$ converges to a point $z_{0}\in Y.$ 
Let $S_{i,j}=f_{g^{n_{i}}\pi (z)}^{n_{j}-n_{i}}$ for 
$1\leq i\leq j.$ 
We set 
$$ V=\{ a\in Y_{\pi (z_{0})}\mid \exists \epsilon >0,\ 
\lim\limits _{i\rightarrow \infty }\sup\limits _{j>i}
\sup\limits _{\xi \in \tilde{B}(a,\epsilon )}
d(S_{i,j}\circ \varphi (\xi ),\ \xi )=0\} ,$$
where $\varphi $ is a map from $Y_{\pi (z_{0})}$ onto 
$Y_{g^{n_{i}}\pi (z)}$ defined by the local parameterization
around $\pi (z_{0}).$  
Then $V$ is a non-empty open proper subset of $Y_{\pi (z_{0})}$, 
and it follows that 
$$ \partial V\subset \hat{J}_{\pi (z_{0})}(f)\cap UH(f).$$ 
\end{thm} 
\begin{cor}
\label{keycor}
Let $(\pi ,Y,X)$ be a $\CCI $-bundle.
Let $f:Y\rightarrow Y$ be a fibered rational map 
over $g:X\rightarrow X.$ 
Assume that $f$ satisfies condition (C1). 
Assume also that, for each $x\in X,\ $ 
the boundary of 
$\hat{J}_{x}(f)\cap UH(f)$ in 
$Y_{x}$ does not separate points 
in $Y_{x}.$ Then,  
for each $z\in Y $ with $z\in F_{\pi (z)},\ $ 
it follows that 
diam$f_{\pi (z)}^{n}(W)\rightarrow 0$
 as $n\rightarrow \infty $, for each 
 small connected neighborhood $W$ of 
 $z$ in $Y_{\pi (z)}$, and that 
 $d(f_{\pi (z)}^{n}(z),UH(f))\rightarrow 0$ 
 as $n\rightarrow \infty .$   
\end{cor}
\begin{proof}
Let $z\in Y$ be a point with 
$z\in F_{\pi (z)}.$ 
Suppose that 
the assumptions of Theorem~\ref{keyth}
are satisfied, and continue with the same notation as that in 
Theorem~\ref{keyth}. 
Let $x_{i}=g^{n_{i}}\pi (z).$ 
Let $D^{i}$ be the 
topological disk in 
$Y_{x_{i}}$ 
satisfying  
condition(C1). 
We may assume that 
there exists a non-empty open set A 
in $\CCI $ such that 
$A\subset i_{x_{i}}^{-1}(D_{i})$ for 
each $i\in \NN .$ 
Then, by the definition of $V$, 
it must be that 
$i_{\pi (z_{0})}(A)\subset Y_{\pi (z_{0})}\setminus V.$ 
This contradicts 
 the assumption of the corollary.  
Hence, 
it must be that diam$f_{\pi (z)}^{n}(W)\rightarrow 0$
 as $n\rightarrow \infty $, for each 
 small connected neighborhood $W$ of 
 $z$ in $Y_{\pi (z)}.$ 
 Then, by Lemma~\ref{backsmall},\ it follows that 
 that  $d(f_{\pi (z)}^{n}(z),UH(f))\rightarrow 0$ 
 as $n\rightarrow \infty .$ 
\end{proof}
\subsection{Proofs of results on fibered rational maps}
\label{Proofs on f}
In this section, we present proofs of the results from 
Section~\ref{Results on f}. 
We begin with the following. \\ 
 {\bf Proof of Theorem~\ref{unifperf1}.}
First, we prove 
that there exists a constant $C>0$ such that 
for each $x\in X$,\ $J_{x}$ is $C$-uniformly perfect.
Since, for each $x\in X,\ $ 
$J_{x}$ is a non-empty perfect set
(Proposition~\ref{potentialprop}),\  
it has uncountably many points. Combined with the 
lower semicontinuity of the map $x\mapsto J_{x}$ 
(statement \ref{potentialprop11}
in Proposition~\ref{potentialprop}) and the 
compactness of $X$,\ this suggests that, 
 for any $x\in X,\ $ one can take four points 
 $z_{x,1},z_{x,2},z_{x,3}$ and $z_{x,4}$ 
 in $J_{x}$ so that 
 $d(z_{x,i},z_{x,j})>C_{1}$, whenever 
 $i\neq j$ and $x\in X$, for some 
 constant $C_{1}$ independent of 
 $(i,j)$ and $x\in X.$ 

 Suppose that there exists a sequence of 
 annuli $\{ D_{j}\} $ with 
 $D_{j}\subset Y_{x_{j}},\ x_{j}\in X $ 
 such that $D_{j}$ separates $J_{x_{j}}$, for each $j$ 
 and mod $D_{j} \rightarrow \infty $ as $j\rightarrow \infty .$ 
Let $D_{j}'$ and $D_{j}''$ be the two components 
of $Y_{x_{j}}\setminus D_{j}.$
We may assume that  
\begin{equation}
\label{unifperfeq}
\mbox{diam }D_{j}''
\rightarrow 0 \mbox{ as }j\rightarrow \infty .
\end{equation}
For, by the existence of $\{ z_{x,i}\} _{x,i},\ $ 
it may be assumed that $\inf_{j\in \NN }\mbox{diam }D_{j}'>0.$ 
Then, since mod $D_{j}\rightarrow \infty $ as $j\rightarrow 
\infty ,\ $ by Lemma 6.1 on p. 34 in \cite{LV} it follows that 
diam $D_{j}''\rightarrow 0$ as $j\rightarrow \infty .$ 

 It may also be assumed that $(\mbox{int}D_{j}'')\cap J_{x_{j}}\neq \emptyset $, 
 for each $j\in \NN $.
Hence, there exists a smallest positive integer    
 $n_{j}$ such that diam $f^{n_{j}+1}(D_{j}'')\geq C_{1}.$ 
 Then, there exists a constant $l_{0}$ 
 such that $l_{0}C_{1}<\mbox{diam }f^{n_{j}}(D_{j}'')$, 
 for each $j.$ Since diam $f^{n_{j}}(D_{j}'')<C_{1},\ $ 
 there exist three distinct points in 
 $\{ z_{x_{j}',i}\} _{i=1,\ldots ,4}$ none of 
 which 
  belongs to $f^{n_{j}}(D_{j}''),\ $ 
 where $x_{j}'=g^{n_{j}}(x_{j}).$ 
 Since $D_{j}\subset F_{x_{j}},\ $ it follows 
 that none of these three points belongs  
 to $f^{n_{j}}(D_{j}),\ $ or to $f^{n_{j}}(D_{j}\cup D_{j}'').$ 
 Let $\varphi _{j}:\{ |z|<1\} \rightarrow D_{j}\cup D_{j}''$ 
 be a 
 Riemann map such that $\varphi _{j}(0)=y_{j}\in 
 D_{j}''$ (note that we may assume that 
 int $ D_{j}'\neq \emptyset 
 $ for each $j$ ). Then, 
 from the above, it follows that, 
 if we set $\alpha _{j}=i_{x_{j}'}^{-1}f^{n_{j}}\varphi _{j}
 :\{ |z|<1\} \rightarrow \CCI $, then 
 $\{ \alpha _{j}\} _{j\in \NN }$ is normal in $\{ |z |<1\} .$ 
 But this causes a contradiction, because 
 diam $\varphi _{j}^{-1}(D_{j}'')\rightarrow 0$ 
 as $j\rightarrow \infty ,\ $ which follows 
 from 
 mod $\varphi _{j}^{-1}(D_{j})=$ mod $D_{j} \rightarrow \infty $ 
 as $j\rightarrow \infty ,\ $ 
 and $l_{0}C_{1}<\mbox{ diam }f^{n_{j}}(D_{j}'')$, 
 for each $j.$  Hence, we have shown 
 that there exists a constant $C>0$ such that 
for each $x\in X$,\ $J_{x}$ is $C$-uniformly perfect.   

 Next, we prove that there exists a constant $C>0$ 
 such that for each $x\in X$,\ $\hat{J}_{x}$ is 
 $C$-uniformly perfect.
 Suppose that there exists a sequence of 
 annuli $\{ D_{j}\} $ with 
 $D_{j}\subset Y_{x_{j}},\ x_{j}\in X $ 
 such that $D_{j}$ separates $\hat{J}_{x_{j}}$, for each $j$, 
 and mod $D_{j} \rightarrow \infty $ as $j\rightarrow \infty .$ 
Let $D_{j}'$ and $D_{j}''$ be the two components 
of $Y_{x_{j}}\setminus D_{j}.$ As in the previous 
paragraph,\ it may be assumed that 
diam $D_{j}''\rightarrow 0$ as $j\rightarrow \infty .$ 

Fix any $j\in \NN .$ Let $y\in D_{j}''\cap \hat{J}_{x_{j}}$ be 
a point. There exists a sequence 
$((x_{j,n},y_{j,n}))_{n}$ in $X\times Y$ with 
$y_{j,n}\in J_{x_{j,n}}$, for each $n\in \NN $, such that 
$(x_{j,n},y_{j,n})\rightarrow (x,y)$ as $n\rightarrow \infty .$ Then 
we may assume that there exists a number $n(j)\in \NN $ such that 
$i_{x_{j},n(j)}^{-1}(J_{x_{j,n(j)}})
\subset i_{x_{j}}^{-1}(D_{j}'\cup D''_{j}).$ 
  Since $i_{x_{j},n(j)}^{-1}(J_{x_{j,n(j)}})\cap i_{x_{j}}^{-1}(D''_{j})
  \neq \emptyset $ (take 
  $n(j)$ sufficiently large),\ by the previous paragraph  
  it must be that $i_{x_{j},n(j)}^{-1}(J_{x_{j,n(j)}})
  \subset i_{x_{j}}^{-1}(D''_{j})$, for large $j.$ 
  However, this contradicts the existence of 
  $\{ z_{x,i}\}_{x,i}$, because    
  it is also the case that diam $D_{j}''\rightarrow 0 $ 
  as $j\rightarrow \infty .$ Hence, we have proved the 
  first and the second statements in \ref{unifperf1-1}.
  The third statement in \ref{unifperf1-1} follows 
  from the uniform perfectness of $J_{x}$ ,\ 
  the continuity of $\omega _{x}$, and Theorem 
  7.2 in \cite{Su}.

 Next, we prove statement \ref{unifperf1-2}.
 Suppose the point $z$ belongs to the boundary 
 of $\hat{J}_{x}$ with respect to the 
 topology of $Y_{x}$, where $x=\pi (z).$ 
 Under a coordinate exchange,\ 
 the map $f_{x}^{n}$ around $z$ is 
 conjugate to $\alpha (z)=z^{l}$ 
 for some $l\in \NN .$ Since 
 $f_{x}^{n}(Y_{x}\setminus \hat{J}_{x})
 \subset Y_{x}\setminus \hat{J}_{x},\ $ 
 there exists an annulus $A$ around $z$ 
 in $Y_{x}$ that separates $\hat{J}_{x}$ 
 and is isomorphic to a round annulus 
 $A'=\{ r<|z|<R\} $ in the above coordinate.
 Then, 
 $\mbox{mod }(f_{x}^{ns}(A))=
 \mbox{mod }(\alpha ^{s}(A'))\rightarrow \infty $ 
 as $s\rightarrow \infty .$ 
 In addition, $f_{x}^{ns}(A)$ separates $\hat{J}_{x}$, 
 for each $s\in \NN .$ This contradicts the fact that 
 $\hat{J}_{x}$ is uniformly perfect.  
 \qed \\ 

 Before proving Theorem~\ref{semihypJohn},
we need the following.
\begin{cor}(Corollary of Theorem~\ref{unifperf1})
\label{unifperfcor}
Let $(\pi,Y=X\times \CCI,\ X)$ be a trivial
$\CCI $-bundle. Let $f:Y\rightarrow Y$ be a fibered
polynomial map over $g:X\rightarrow X$ such that
$d(x)\geq 2$,\ for each
$x\in X$. Let $R>0$ be a number
such that for each $x\in X$,
$\pi _{\CCI }(J_{x}(f))\subset
\{ z\mid d(z,\infty )>R\}$,
where $\pi _{\CCI }:Y\rightarrow \CCI $ denotes
the projection and $d$ denotes the
spherical distance. (Note that such an $R$ exists,
since $d(x)\geq 2$ for each $x\in X.$)
Let $\rho _{x}(z)|dz|$ be the hyperbolic metric
on $\pi _{\CCI }(A_{x}(f))$. Let $\delta _{x}(z)=
\inf _{w\in \partial (\pi _{\CCI }(A_{x}(f)))}
|z-w|$ for each $z\in \pi _{\CCI }(A_{x}(f))\cap \CC.$
Then, there exists a positive constant $C$ depending only on
$f$ and $R$ such that for each $x\in X$ and
each $z\in \pi _{\CCI }(A_{x}(f))\cap \{ z\mid
d(z,\infty )>R\} $,\ we have
$C\leq \rho _{x}(z)\delta _{x}(z)\leq 1$.

\end{cor}
\begin{proof}
$\rho _{x}(z)\delta _{x}(z)\leq 1$ follows easily from
the Schwarz lemma.

Next, as in the proof of Theorem~\ref{unifperf1},
for any $x\in X$, we can take two points
$z_{x,1}$ and $z_{x,2}$ in $J_{x}(f)$ so that
$d(z_{x,1},z_{x,2})>c_{0}$, whenever
$x\in X$, for some positive constant $c_{0}$
independent of $x\in X.$
Let $\psi _{x}(z)$ be a M\"{o}bius transformation
such that $\psi _{x}(z_{x,1})=\infty $ and $\psi _{x}$
preserves the spherical metric.
Let $B_{x}=\psi _{x}(\pi _{\CCI }(A_{x}(f)))\ (\subset \CC ).$
By Theorem~\ref{unifperf1} and Theorem 2.16 in \cite{Su},
there exists a positive constant
$c_{1}$ such that for each $x\in X$ and
each $z\in B_{x}$,
$\rho _{x,1}(z)\inf _{w\in \partial B_{x}}|z-w|\geq c_{1}$,
where $\rho _{x,1}(z)|dz|$ denotes the hyperbolic metric on
$B_{x}$.

Let $z\in \pi _{\CCI }(A_{x}(f))\cap \{ y\mid d(y,z_{x,1})>c_{0}/2,\ 
d(y,\infty )>R\}$. Then, we have
$\rho _{x}(z)=\rho _{x,1}(\psi _{x}(z))|\psi _{x}'(z)|$.
Since $\psi _{x}(z)\in \{ y\mid  
d(y,\infty )>c_{0}/2,\ d(y,\psi _{x}(\infty ))>R\}$ and
$\psi _{x}^{-1}(\{ y\mid d(y,\psi _{x}(\infty ))>R\})
\subset \{y\mid d(y,\infty )>R\}$,
by the Cauchy formula, we have
$|\psi _{x}'(z)|=|(\psi _{x}^{-1})'(\psi _{x}(z))|^{-1}\geq c_{2}$,
where $c_{2}$ is a positive constant independent of
$z$ and $x$.

Next, let $w_{0}\in \partial (\pi _{\CCI }(A_{x}(f)))
=\pi _{\CCI }(J_{x}(f))$ be a point
such that $\delta _{x}(z)=|z-w_{0}|$.

Suppose case (1): $w_{0}\in \{y\mid
d(y,z_{x,1})\leq c_{0}/4\}.$ Then,
$|z-w_{0}|\geq c_{3}$, where
$c_{3}$ is a positive constant depending only on
$c_{0}$. Further,
$\inf _{w\in \partial B_{x}}|\psi _{x}(z)-w|
\leq |\psi _{x}(z)-\psi _{x}(z_{x,2})|\leq c_{4}$,
where $c_{4}$ is a positive constant depending only on $c_{0}$.
Hence, $\delta _{x}(z)\geq \frac{c_{3}}{c_{4}}\inf _{w\in \partial
B_{x}}|\psi _{x}(z)-w|$.

Suppose case (2): $w_{0}\in \{y\mid
d(y,z_{x,1})>c_{0}/4\}$. Let $\gamma$ be the
Euclidean segment connecting $z$ and $w_{0}$. Then,
since $|w_{0}-z|=\delta _{x}(z)$,\ 
we obtain $d(z_{x,1},\gamma )\geq c_{5}$, where $c_{5}$ is a
positive constant depending only on $c_{0}$,
which implies $d(\infty,\psi _{x}(\gamma ))$ $\geq c_{5}$.
Hence, by the Cauchy formula, we have
$\inf _{w\in \partial B_{x}}|\psi _{x}(z)-w|
\leq |\psi _{x}(z)-\psi _{x}(w_{0})|
\leq \sup _{w\in \gamma}|\psi _{x}'(w)|\cdot |z-w_{0}|
\leq c_{6}\delta _{x}(z)$,
where $c_{6}$ is a positive constant independent
of $z$ and $x_{0}$.

From these arguments, we find that there exists a
positive constant $c_{7}$ depending only on $f$ and $R$, such that
for each $z\in \pi _{\CCI }(A_{x}(f))\cap \{y\mid
d(y,z_{x,1})>c_{0}/2,\ d(y,\infty )>R\}$, we have
$\rho _{x}(z)\delta _{x}(z)\geq c_{7}$.
Similarly, we find that there exists a positive constant
$c_{8}$ depending only on $f$ and $R$, such that
for each $z\in \pi _{\CCI }(A_{x}(f))\cap
\{y\mid d(y,z_{x,2})>c_{0}/2,\ d(y,\infty)>R\}$,
we have $\rho _{x}(z)\delta _{x}(z)\geq c_{8}.$
Since $d(z_{x,1},z_{x,2})>c_{0}$, we obtain the
statement of the corollary.

\end{proof}
{\bf Proof of Theorem~\ref{semihypJohn}.}
This statement can be proved using an approach similar to that in
\cite{CJY}, using Theorem 2.14 (3) in \cite{S4},
Corollary~\ref{unifperfcor},
Lemma~\ref{ylem1} and Lemma~\ref{backsmall}.

The procedure is as follows: first,
for any $(x,y)\in X\times \CC $, let
$G_{x}((x,y)):=\lim_{n\rightarrow \infty }\frac{1}{d_{n}(x)}
\log ^{+}|\pi _{\CCI }(f_{x}^{n}((x,y)))|$, 
where $\log ^{+}(a):=\max \{ \log a,0\}$
for $a\in \RR$ with $a>0.$
Then, by the arguments in \cite{Se3} (or \cite{Se1}),
on $A_{x}(f)$, $G_{x}$ equals the
Green's function with pole $(x,\infty )$,
$G_{g(x)}(f_{x}(z))=d(x)G_{x}(z)$ for
any $x\in X$ and $z\in Y_{x}$,
$(x,y)\mapsto G_{x}((x,y))$ is continuous
in $X\times \CC $, and
there exist numbers $r>0$ and $M>0$ such that
for any $(x,y)\in X\times \CC $ with
$|y|>r$, $|G_{x}((x,y))-\log |y||\leq M$.
Further, for any $x\in X$, we can define
a locally conformal map $\psi _{x}$ that is defined on
$\{(x,y)\mid |y|>R\}$, where $R$ is a number independent
of $x\in X$, such that
$\psi _{g(x)}f_{x}=(\psi _{x})^{d(x)}$
and $\log |\psi _{x}|=G_{x}$ on $\{(x,y)\mid
|y|>R\}$. In particular, $G_{x}$ has no
critical points in $\{(x,y)\mid |y|>R\}$.

Next, for each $x\in X$, let $\rho _{x}(z)|dz|$ be the
hyperbolic metric on $A_{x}(f)$. Similarly,
let $\delta _{x}(z)|dz|$ be the quasihyperbolic
metric on $A_{x}(f)\
(\delta _{x}(z):=\inf _{w\in J_{x}}
|\pi _{\CCI }(z)-\pi _{\CCI }(w)|$).
By Theorem 2.14(3) in \cite{S4},
we obtain the following claim:\\
Claim 1: There exist constants $C>0$ and $\alpha >0$
such that for each $x\in X$ and each $z\in A_{x}(f)\cap
(X\times \CC )$ with $\delta _{x}(z)\leq 1$,
we have $\delta _{x}(z)\leq CG_{x}(z)^{\alpha }$.

Next, for each $x\in X$ and each $z\in A_{x}\cap (X\times
\CC )$, we take a Green's line (a curve
of steepest descent for $G_{x}$)
$\gamma _{z}$
running from $z$
to ``$J_{x}(f)$", such that
$f_{x}(\gamma _{z})= \gamma _{f_{x}(z)}$.
(We choose a direction (left or right) consistently
at every critical point encountered.)
Then, we show the following:\\
Claim 2: Let $0<C_{1}<C_{2}<\infty,\ 
0<C_{3}<C_{4}<\infty $ be given numbers.
Let $\tau :[a,b]\rightarrow A_{x}(f)\ (0\leq a<b<\infty )$
be a curve that parameterizes part of a Green's line,
such that $\tau ([a,b])\subset \{ z\in A_{x}(f)\mid
C_{1}<\delta _{x}(z)<C_{2}\}$ and
$C_{3}\leq l_{E}\tau <C_{4}$, where
$l_{E}$ denotes the Euclidean length.
Then, there exist constants $C_{5}>0$ and $0<\lambda _{1}<1$,
which depend only on $f, C_{1}, C_{2}, C_{3}$, and $C_{4}$,
such that $|\tau (a)-\tau (b)|\geq C_{5}$ and
$G_{x}(\tau (b))<\lambda _{1}G_{x}(\tau (a)).$

The proof of this claim is immediate, since we assume that
either $J_{x}(f)$ is connected for each $x\in X$
(then $f_{x}$ and $G_{x}$ have no critical points in $A_{x}(f)$)
or condition \ref{semihypJohnh2} in the 
assumption of Theorem~\ref{semihypJohn} 
is satisfied. 
Using Claim 2,\ 
$\inf \{ d(a,b)\mid a\in \pi _{\CCI }(J_{x}),\ 
b\in \pi _{\CCI }(A_{x}(f)\cap C(f)),\ x\in X\} >0$,\ 
 the Koebe distortion theorem, and
Corollary~\ref{unifperfcor}, we easily obtain
the following:\\
Claim 3: Let $0<C_{6}<C_{7}<\infty $ and
$0<C_{8}<\infty $ be given numbers.
Let $\tau :[a,b]\rightarrow A_{x}(f)$ be a
curve that parameterizes part of a Green's line.
Suppose $\tau ([a,b])\subset
\{ z\in A_{x}(f)\mid |\pi _{\CCI }(z)|\leq C_{8}\}$
and $C_{6}\leq l_{\rho _{x}}\tau \leq C_{7}$
($l_{\rho _{x}}$ denotes the $\rho _{x}$-length).
Then, there exist constants $C_{9}>0$
and $0<\lambda _{2}<1$ that depend only on
$f, C_{6}, C_{7}$, and $C_{8}$, such that
$C_{9}\leq d_{\rho _{x}}(\tau (a),\tau (b))\leq C_{7}$
($d_{\rho _{x}}$ denotes the distance induced by
$\rho _{x}(z)|dz|$) and
$G_{x}(\tau (b))<\lambda _{2}G_{x}(\tau (a)).$

Next, we will show the following:\\
Claim 4: There exist constants $C_{10}>0, C_{11}>0$, and
$\alpha >0$ such that for any $x\in X$ and
any $z\in A_{x}(f)$ with $G_{x}(z)\leq C_{10}$,
we have $l_{E}\gamma _{z}\leq C_{11}G_{x}(z)^{\alpha }$.
In particular, any Green's line $\gamma _{z}$ lands at
a point in $J_{x}(f).$

To show this claim, take a number $C_{10}$ such that
for any $x\in X$,
$\{z\in A_{x}(f)\mid G_{x}(z)\leq C_{10}\}
\subset \{z\in A_{x}(f)\mid \delta _{x}(z)\leq 1\}$.
Next, for any two points $w_{1}$ and $w_{2}$ in a
Green's line $\gamma _{z_{0}}\ (z_{0}\in A_{x}(f)) $ with
$G_{x}(w_{2})\leq G_{x}(w_{1})$,
let $\gamma _{z_{0}}(w_{1},w_{2})$ be the arc from
$w_{1}$ to $w_{2}$ along $\gamma _{z_{0}}$.
Let $\{ z_{n}\} $ be a sequence in $\gamma _{z_{0}}$
such that
$l_{\rho _{x}}\gamma _{z_{0}}(z_{n},z_{n+1})=c$ for each
$n\in \NN$,
where $c$ is a positive constant that is sufficiently
small. Then, by Corollary~\ref{unifperfcor}, Claim 1,
and Claim 3, we obtain
$l_{E}\gamma _{z_{0}}=\sum _{n}
l_{E}\gamma _{z_{0}}(z_{n},z_{n+1})\leq
\mbox{Const.}\sum _{n}
\int _{\gamma _{z_{0}}(z_{n},z_{n+1})}
\delta _{x}(z)\rho _{x}(z)|dz|
\leq \mbox{Const.}\sum _{n}
\int _{\gamma _{z_{0}}(z_{n},z_{n+1})}
\delta _{x}(z_{n})\rho _{x}(z)|dz|
\leq \mbox{Const.}\sum _{n}\delta _{x}(z_{n})
\leq \mbox{Const.}\sum _{n}G_{x}(z_{n})^{\alpha }$
$\leq \mbox{Const.}\sum _{n}\lambda ^{n\alpha }G_{x}(z_{0})^{\alpha}$
$\leq \mbox{Const.}G_{x}(z_{0})^{\alpha }$,
where, $\alpha >0$ and $0<\lambda <1$ are constants and
Const. denotes a constant. Hence, we obtain the claim.

Next, we show the following:\\
Claim 5: There exist constants $\delta '>0$ and $p\in \NN$
such that for any $x\in X$ and any $(w,z)\in A_{x}(f)\times
A_{x}(f)$ with $w\in \gamma _{z},\ 
l_{\rho _{x}}\gamma _{z}(z,w)\geq p$ and
$G_{x}(z)<\delta '$,
we have $\delta _{x}(w)\leq \frac{1}{2}\delta _{x}(z)$.

This claim can be shown using the argument on page 11 in
\cite{CJY}, using Claims 3 and 4, and the
distortion lemma for proper holomorphic maps.
Next, we show the following:\\
Claim 6: There exist constants $\delta '>0$ and $C_{12}>0$
such that for any $x\in X$ and any $(w,z)\in A_{x}(f)\times
A_{x}(f)$ with $w\in \gamma _{z}$ and $G_{x}(z)\leq \delta '$,
we have
$\delta _{x}(z)\geq C_{12}|\pi _{\CCI}(z)-\pi _{\CCI }(w)|$.

This claim can be shown using Claim 5 and Corollary~\ref{unifperfcor}.
Finally, we show the following:\\
Claim 7: There exists a constant $C>0$ such that
for any $x\in X$ and any $(w,z)\in A_{x}(f)\times
A_{x}(f)$ with $w\in \gamma _{z}$, we have
$\delta _{x}(z)\geq C|\pi _{\CCI }(z)-\pi _{\CCI }(w)|$.

This claim can easily be shown using Claim 6 and
the fact that there exist numbers $r>0$ and $M>0$ such that
for any $(x,y)\in X\times \CC $ with
$|y|>r$, $|G_{x}((x,y))-\log |y||\leq M$.
Hence, we have shown Theorem~\ref{semihypJohn}.
\qed \\

\noindent {\bf Proof of Theorem~\ref{porositythm}.}
For any $y'\in \bigcup _{x\in X}J_{x}$
and $r>0,\ $ we set 
$$ h(y',r)=\sup \{ s\mid \exists y''\in Y_{\pi (y')},\ 
\tilde{B}(y'',s)\subset F_{\pi (y')}\cap \tilde{B}(y',r)\} $$ 
and $h(r)=\inf \{ h(y',r)\mid y'\in \bigcup _{x\in X}J_{x}\} .$
By Theorem 2.14 and Remark 6 in \cite{S4},\  it follows that 
$\tilde{J}(f)=\bigcup _{x\in X}J_{x}.$ 
Furthermore, by condition (C1), int $J_{x}=\emptyset $, 
for any $x\in X.$ Hence, it follows that
$h(r)>0$ for any $r>0.$ 

 Since $f$ is semi-hyperbolic and 
 satisfies condition (C1),\ by Lemma~\ref{keythmlem}  
there exists a positive number 
$\delta _{1}$ and a number $N\in \NN $ 
such that, for any $y'\in \tilde{J}(f),\ 
0<\delta \leq \delta _{1},\ $ $n\in \NN $, 
any $V\in c(\tilde{B}(y',2\delta ),\ f^{n}),\ $ 
we have that $V$ is simply connected and that 
$\deg (f^{n}:V\rightarrow \tilde{B}(y' ,2\delta ))\leq N.$ 

 Let $y\in \tilde{J}(f)$ and $r>0.$ We set 
 $B_{n}=f^{n}(\tilde{B}(y,r))$ and 
 $y_{n}=f^{n}(y)$, for each $n\in \NN .$ 
 Since $y\in J_{\pi (y)},\ $ 
 there exists a smallest positive integer 
 $n_{0}$ such that 
 diam $B_{n_{0}+1}>\delta _{1}.$ Then there exists a 
 constant $l_{0}$ such that 
 $l_{0}\delta _{1}<$diam $B_{n_{0}}.$ 
  By Corollary 2.3 in \cite{Y},\ 
  there exists a constant $K$ that depends only on 
  $N$ and a ball 
  $\tilde{B}(y_{n_{0}},r_{0})\subset B_{n_{0}}$, with 
 $r_{0}\geq $diam $B_{n_{0}}/K\geq \frac{l_{0}\delta _{1}}{K}$ 
such that the element 
$W \in c((\tilde{B}(y_{n_{0}},r_{0}),\ f^{n_{0}})$
containing $y$ is a subset of $\tilde{B}(y,r).$ 
There exists a ball $\tilde{B}(y',\frac{2}{3}h(r_{0}))$ 
included in $\tilde{B}(y_{n_{0}},r_{0})\cap F_{\pi (y_{n_{0}})}.$ 

 Let $D_{0}$ be an element of  
 $c((\tilde{B}(y',\frac{1}{2}h(r_{0}))),\ f^{n_{0}})$
 contained in $W.$ It follows that 
 $D_{0}\subset F_{\pi (y)}.$ 
 Let $y''\in D_{0}\cap (f^{n_{0}})^{-1}\{ y'\} $ 
 be a point. Then, by Corollaries 1.8 and 1.9 in \cite{S4},\ 
 it follows that Dist $(\partial D_{0},y'')\leq M$, for some $M$ 
 that depends only on $N.$ 
 Further, by the method used in Corollary 1.9 in \cite{S4}
 it follows that 
 $\min _{z\in \partial D_{0}} d_{\pi (y'')}(y'',z)\asymp 
  \max _{z\in \partial E}d_{\pi (y'')}(y'',z)$,
  where $E\in c(\tilde{B}(y_{n_{0}},\delta _{1}),\ f^{n_{0}})$ 
  with $y'' \in E.$ Since $\tilde{B}(y,r)\subset E$, it follows that 
 diam $D_{0}\asymp r.$ Combining this with Dist $(\partial D_{0},y'')\leq M$
 leads to the conclusion that there exists a constant $0<k<1$ that does not 
 depend on $y$ and $r$ such that 
 $B(y'',kr)\subset D_{0}\subset F_{\pi (y)}.$ 
\qed \\   

To show Theorem~\ref{bundlemeas},\ we present 
some notation and propositions.
\begin{df}({\bf Conical set for fibered rational maps})
Let $(\pi , Y,X)$ be a $\CCI $-bundle. 
Let $f:Y\rightarrow Y$ be a fibered rational map 
over $g:X\rightarrow X.$ 
Let $N\in \NN $ and $r>0.$ 
We denote by $\tilde{J}_{con}(f,N,r)$ 
the set of points $z\in \tilde{J}(f)$ 
such that, 
for any $\epsilon >0,\ $ there exists a 
positive integer $n$ such that 
 the 
 element 
 $U\in c(\tilde{B}(f^{n}(z),r),\ 
 f^{n}|_{Y_{\pi (z)}})$ containing $z$ satisfies the 
 following conditions:
 \begin{enumerate}
 \item diam $U\leq \epsilon, $
 \item $U$ is simply connected,\ and 
 \item $\deg (f^{n}:U\rightarrow  
 \tilde{B}(f^{n}(z),r))\leq N.$
 \end{enumerate} 
We set $\tilde{J}_{con}(f,N)=\bigcup _{r>0}
\tilde{J}_{con}(f,N,r)$ and 
$\tilde{J}_{con}(f)=\bigcup _{N\in \NN }
\tilde{J}_{con}(f,N).$  
\end{df}
\begin{prop}
\label{measprop1}
Let $(\pi , Y,X)$ be a $\CCI $-bundle.
Let $f:Y\rightarrow Y$ be a 
fibered rational map over $g:X\rightarrow X.$ 
Suppose that $\hat{J}_{x}$ has no interior points with respect to 
the topology in $Y_{x}$, 
for each $x\in X.$ Then the two-dimensional Lebesgue 
measure of $\tilde{J}_{con}(f)\cap J_{x}$ is 
equal to zero, for each $x\in X.$ 
\end{prop}
\begin{proof}  
Fix $N\in \NN .$ Suppose that there exists a point $x\in X$ such that 
$\tilde{J}_{con}(f,N)\cap J_{x}$
has positive measure. Then there exists a 
Lebesgue density point $y\in \tilde{J}_{con}(f,N)\cap J_{x}.$
Let $y_{m}=f_{x}^{m}(y)$ and $x_{m}=g^{m}(x)$, 
for any $m\in \NN .$  Let $\delta >0$ be a number 
such that $y\in \tilde{J}_{con}(f,N,\delta ).$
 Let $U_{m},U_{m}'$ be the respective
elements of $c(\tilde{B}(y_{m},\delta /2 ),\ f_{x}^{m}),\ $
$c(\tilde{B}(y_{m},\delta ),\ f_{x}^{m})$  
containing $y$. Since 
$y\in \tilde{J}_{con}(f,N,\delta ),\ $  
there exists a
 subsequence $(n)$ in $\NN $ with $n\rightarrow \infty $ 
  such that $U_{n}'$ is simply connected,\ 
$\deg (f_{x}^{n}:U_{n}'\rightarrow  
\tilde{B}(y_{n},\delta ))\leq N $, for each $n$, and 
diam $U_{n}'\rightarrow 0$ as $n\rightarrow \infty .$  
By Corollary 1.9 in \cite{S4},  for any 
local parameterization $i_{x},\ $ 
\begin{equation}
\label{bundlemeaseq0}
 \lim\limits _{n\rightarrow \infty } \frac{m(i_{x}^{-1}(U_{n}\cap J_{x}))}
{m(i_{x}^{-1}(U_{n}))}=1,\ 
\end{equation} where 
$m$ denotes the spherical measure of $\CCI .$  
Using an argument from the proof of Theorem 4.4 
in \cite{S4},\ from (\ref{bundlemeaseq0}) we 
can show that 
\begin{equation}
\label{bundlemeaseq1}
\lim\limits _{n\rightarrow \infty } \frac{m(\i_{x_{n}}^{-1}
(\tilde{B}(y_{n},\delta /2 )\cap F_{x_{n}}))}
{m(i_{x_{n}}^{-1}(\tilde{B}(y_{n}, \delta /2)))}=0, 
\end{equation} 
where $i_{x_{n}}$ denotes a local parameterization. 
There exists a subsequence $(n_{j})$ of $(n),\ $ 
a point $y_{\infty }\in Y$, and a point $x_{\infty }\in X$ 
such that $y_{n_{j}}\rightarrow y_{\infty }$ and 
$x_{n_{j}}\rightarrow x_{\infty }$ as 
$j\rightarrow \infty .$ By (\ref{bundlemeaseq1}),
it follows  
that $\tilde{B}(y_{\infty },\ \delta /2)\subset 
\hat{J}_{x_{\infty }}.$ On the other hand,\ 
by assumption, for any 
$a\in X,\ $ $\hat{J}_{a}$ has no interior point. 
This is a contradiction.   

\end{proof}

\begin{prop}
\label{measprop2}
Let $(\pi , Y,X)$ be a $\CCI $-bundle.
Let $f:Y\rightarrow Y$ be a 
fibered rational map over $g:X\rightarrow X.$ 
Suppose that $f$ satisfies condition (C1). Then 
the following hold.
\begin{enumerate}
\item 
$$\tilde{J}_{good}(f)\cap \bigcup _{x\in X}J_{x}\subset \tilde{J}_{con}(f).$$ 
\item 
If it is assumed, furthermore, that, for each $x\in X,\ $ 
the boundary of $\hat{J}_{x}(f)\cap UH(f)$ 
in $Y_{x}$ does not separate points in $Y_{x},\ $ 
then 
$$ \tilde{J}_{good}(f)\subset \tilde{J}_{con}(f)\cap 
\bigcup _{x\in X}J_{x}.$$ 
\end{enumerate}
\end{prop}
\begin{proof} We will first prove the first statement. 
Let $z\in \bigcup _{x\in X}J_{x}$ be a 
point such that 
$\limsup _{n\rightarrow \infty }d(f^{n}(z),UH(f))>0.$
For each $m\in \NN $, let $z_{m}=f^{m}(z)$ 
and $x_{m}=\pi (f^{m}(z)).$ For each $m\in \NN $ 
and 
each $r>0$, let $U_{m}(r),U_{m}'(r)$ 
respectively, be the elements of $c(\tilde{B}(z_{m},r/2),\ f_{\pi (z)}^{m}), 
c(\tilde{B}(z_{m},r),\ f_{\pi (z)}^{m})$ containing 
$z$.
 There exists a positive number $\delta ,\ $ 
a positive integer $N$, and 
 a sequence $(n)$ in $\NN $ 
 such that 
 $\deg (f_{\pi (z)}^{n}:U_{n}'(\delta )\rightarrow 
 \tilde{B}(z_{n},\delta ))\leq N.$ 
 By Lemma~\ref{keythmlem},\ taking $\delta $ sufficiently small 
 it can be assumed that $U_{n}'(\delta ) $ is simply connected. 

  Suppose that diam $ (U_{n}(\delta ))$ does not 
  tend to zero as $n\rightarrow \infty $ 
  in $(n). $ Then, 
  by the distortion lemma for proper 
  maps, there exists a subsequence 
  $(n_{j})$ of $(n)$  with $n_{j}\rightarrow \infty $ 
  and a positive number $r$ such that 
  $U_{n_{j}}(\delta )\supset \tilde{B}(z,r) $, 
  for each $j.$ Hence, 
\begin{equation}
\label{measprop2eq1}
f^{n_{j}}(\tilde{B}(z,r))\subset \tilde{B}(f_{n_{j}}(z),\delta )
\end{equation}  
for each $j.$ By condition (C1),\ if 
$\delta $ is sufficiently small, then (\ref{measprop2eq1}) 
contradicts that 
$z\in \cup _{x\in X}J_{x}.$ Hence, 
it follows
that diam $U_{n}(\delta )\rightarrow 0$ 
as $n\rightarrow \infty $ in $(n).$ 
Therefore, $z\in \tilde{J}_{con}(f).$  

 The second statement follows from Corollary~\ref{keycor} and 
 the first statement.   
\end{proof}
\begin{cor}
\label{meascor1}
Let $(\pi , Y,X)$ be a $\CCI $-bundle.
Let $f:Y\rightarrow Y$ be a 
fibered rational map over $g:X\rightarrow X.$ 
Suppose that 
$\tilde{J}(f)=\bigcup _{x\in X}J_{x}$, and that
 $f$ satisfies condition (C1).
Then, for each $x\in X$, it follows that 
$d(f_{x}^{n}(y), UH(f))\rightarrow 0,\ \mbox{ as }
n\rightarrow \infty ,$ 
for almost every $y\in J_{x}$, with respect to the 
Lebesgue measure in $Y_{x}.$ 
\end{cor}
\begin{proof}
By condition (C1), $\hat{J}_{x}=J_{x}$ has no interior 
points, for each $x\in X.$ By Proposition~\ref{measprop1} 
and Proposition~\ref{measprop2},\ 
the following statement results.   
\end{proof}
{\bf Proof of Theorem~\ref{bundlemeas}.}
Suppose that there exists a point 
$z\in \tilde{J}(f)$ such that 
$z\in F_{\pi (z)}.$ By Corollary~\ref{keycor},\ 
for each small connected neighborhood 
$W$ of $z$ in $F_{\pi (z)}$, it follows that 
diam $f^{n}(W)\rightarrow 0$ 
and $d(f^{n}(z),UH(f))\rightarrow 0$ as 
$n\rightarrow \infty .$ However, 
by conditions three and four in the assumptions 
of our theorem,\ this causes a contradiction. 
Hence, we have shown that $\tilde{J}(f)=\bigcup _{x\in X}J_{x}.$ 
By Corollary~\ref{meascor1}, it also follows that 
the two-dimensional Lebesgue measure 
of $J_{x}\setminus \bigcup _{n\in \NN }f^{-n}(UH(f))$ 
is equal to zero.   
\qed \\ 

\noindent  {\bf Proof of Theorem~\ref{rigidity}.}
Let $x'=v(x).$ We now show the following.

 {\bf Claim:} There exists a constant $C>0$ such that, 
 for each $y\in \CCI ,\ $ 
 \begin{equation}
 \label{rigidityeq}
  \limsup\limits _{r\rightarrow 0}
\sup \{ \frac{d_{x'}(u_{x}(i_{x}(y)),\ u_{x}(i_{x}(y')))}
        { d_{x'}(u_{x}(i_{x}(y)),\ u_{x}(i_{x}(y'')))}
  \mid d(y,y')=d(y,y'')=r\} \leq C, 
  \end{equation} 
 where $d_{x'}$  is the metric on 
 $\tilde{Y}_{x'}$
 induced by the form $\tilde{\omega }_{x'}$.

 If $i_{x}(y)\in F_{x}(f),\ $ then (\ref{rigidityeq}) holds 
 with a constant $C=C(K)$ that depends only on $K.$ 
 Let $i_{x}(y)\in J_{x}(f)$ be a point. Let $\delta _{1}$ and 
 $N\in \NN $ be the numbers from the proof of Theorem~\ref{porositythm}.
 Let $r>0$ be a number.  We set 
 $B_{n}=f^{n}(i_{x}(B(y,r)))$ and 
 $y_{n}=f^{n}(i_{x}(y))$, for each $n\in \NN .$ 
 Since $i_{x}(y)\in J_{x}(f),\ $ 
 there exists a smallest positive integer 
 $n_{0}$ such that 
 diam $B_{n_{0}+1}>\delta _{1}.$ Furthermore, there exists a 
 constant $l_{0}$ such that 
 $l_{0}\delta _{1}<$diam $B_{n_{0}}.$ 
  By Corollary 2.3 in \cite{Y},\ 
  there exists a constant $A$ that depends only on 
  $N$, and a ball 
  $\tilde{B}(y_{n_{0}},r_{0})\subset B_{n_{0}}$ with 
 $r_{0}\geq $diam $B_{n_{0}}/A\geq \frac{l_{0}\delta _{1}}{A}$ 
such that the component of 
$(f^{n_{0}}i_{x})^{-1}(\tilde{B}(y_{n_{0}},r_{0}))$
containing $y$ is a subset of $B(y,r).$ 
Then there exist constants $\delta '>0$ and $0<c<1$ such that 
$$ \tilde{B}(u(y_{n_{0}}),\ c\delta ' )\subset u(\tilde{B}(y_{n_{0}},r_{0}))
\subset u(B_{n_{0}})\subset \tilde{B}(u(y_{n_{0}}),\ \delta ' ).$$ 
Since $f$ is semi-hyperbolic and satisfies condition(C1),\ the same thing 
holds for $\tilde{f}$. Applying Lemma~\ref{keythmlem},\  
Lemma~\ref{ylem1}, and the Koebe distortion theorem
to the map $\tilde{f}^{n_{0}},\ $ which maps a neighborhood 
$u(i_{x}(B(y,r)))$ of $u(i_{x}(y))$ in $\tilde{Y}_{x'}$ to 
a neighborhood of $u(y_{n_{0}})$ in 
$\tilde{Y}_{\tilde{g}^{n_{0}}(x')},\ $ 
 it follows that (\ref{rigidityeq}), for 
some constant $C'$ that does not depend on $i_{x}(y)\in J_{x}(f).$ Hence, the claim holds.

 By this claim,\ $u_{x}:Y_{x}\rightarrow \tilde{Y}_{x'}$ is quasiconformal. 
 Since the two-dimensional Lebesgue measure of $J_{x}(f)$ is zero,\ 
 which is a consequence of Theorem~\ref{bundlemeas},\ 
 it is $K$-quasiconformal on $Y_{x}.$ 
 Since $\tilde{f}_{v(w)}\circ u_{w}=u_{g(w)}\circ f_{w}$ on 
 $Y_{w}$, 
 and $f_{w}$ and $\tilde{f}_{v(w)}$ are holomorphic for each 
 $w\in X,\ $ it follows that 
 $u_{a'}:Y_{a'}\rightarrow \tilde{Y}_{v(a')}$ is $K$-quasiconformal
 for each $a'\in \cup _{n\in \Bbb{Z}} g^{n}(\{ x\} ).$ 
 Hence, for each $a\in \overline{\cup _{n\in \Bbb{Z}}g^{n}(\{ x\} )}$, 
 the map $u_{a}:Y_{a}\rightarrow \tilde{Y}_{v(a)}$ is also
 $K$-quasiconformal,\ with the same dilatation constant $K.$   
\qed     

\subsection{Proofs of results on rational semigroups}
\label{Proofs on r}
In this section, we prove the results given in 
Section~\ref{Results on r}.
 We now prove Theorem~\ref{unifperrs}.\\ 
{\bf Proof of Theorem~\ref{unifperrs}.}
Let $f:\Lambda ^{\NN }\times \CCI 
\rightarrow \Lambda ^{\NN }\times \CCI $ 
be the fibered rational map over the 
shift map $\sigma :\Lambda ^{\NN }\rightarrow 
\Lambda ^{\NN }  \ (\sigma (h_{1},h_{2},h_{3},\ldots 
)=(h_{2},h_{3},\ldots ) )$
 defined as $f((h_{1},h_{2},\ldots ) ,\ y)=
 ((h_{2},h_{3},\ldots ),\ h_{1}(y)).$ 
 Then, by Theorem~\ref{unifperf1},\ 
 there exist positive constants 
 $C_{0}$ and $C_{1}$ such that 
 each $g\in G$ satisfies the conditions that 
 $J(g)$ is $C_{0}$-uniformly perfect and 
 that diam $J(g)>C_{1}.$ 

 Let $H$ be any subsemigroup of $G.$
 Let $A$ be an annulus that separates 
 $J(H).$ Let $V_{1}$ and $V_{2}$ be 
 two connected components of $\CCI \setminus 
 A.$ Since 
 $J(G)=\overline{\bigcup _{g\in G}J(g)}$ 
 (Corollary 3.1 in \cite{HM1}) ,\ 
 there exist elements $g_{1}$ and $g_{2}$ in 
 $H$ such that 
 $J(g_{1})\cap V_{1}\neq \emptyset $ and 
 $J(g_{2})\cap V_{2}\neq \emptyset .$ 

 If $J(g_{1})\cap V_{2}\neq \emptyset $
 or $J(g_{2})\cap V_{1}\neq \emptyset ,\ $ 
 then mod $A\leq C_{0}.$ 
 If $J(g_{1})\cap V_{2}=\emptyset $ and 
 $J(g_{2})\cap V_{1}=\emptyset ,\ $ 
 then, by Lemma 6.1 on p. 34 in \cite{LV},
 there exists a constant $C=C(C_{1})$ that 
 depends only on $C_{1} $ such that 
 mod $A\leq C.$ The second statement follows 
 from the uniform perfectness of $J(H)$ and 
 Theorem 4.1 in \cite{HM2}.      
\qed \\ 

 To show Theorem~\ref{rsporosity},\ 
we present some notation and lemmas.

\noindent {\bf Notation:} Throughout the rest of this 
section,\ for a fixed generator 
system $\{ h_{1},\ldots ,h_{m}\} $, 
let $f:\Sigma _{m}\times \CCI 
\rightarrow \Sigma _{m}\times \CCI $ be 
 the fibered rational map 
 over the shift map
 $\sigma :\Sigma _{m}\rightarrow \Sigma _{m},\  $
 where $\Sigma _{m}=\{ 1,\ldots ,m\}^{\NN },\ $ 
 associated with the generator system 
 $\{ h_{1},\ldots ,h_{m}\} .$ 
 We set $q_{x}^{(n)}(y)=\pi _{\CCI }(f_{x}^{n}(y))$, 
 for any $(x,y)\in \Sigma _{m}\times \CCI .$ 
  
\begin{lem}
\label{rsporolem0}
Let $E$ be a finite subset of $\CCI .$ 
Let $\langle h_{1},\ldots ,h_{m}\rangle $ 
be a rational semigroup.
Then, for any number $M>0$, there exists a 
positive integer $n_{0}$ such that, 
for any $(n,x,y)\in \NN \times \Sigma _{m}\times E  $ with 
$n\geq n_{0} $ that  
satisfies all of the following conditions:
\begin{enumerate}

\item $q_{x}^{(j)}(y)\in E$ for $j=0,\ldots ,n$ 
  
\item $(q_{x}^{(n)})'(y)\neq 0 $ and 
\item for any $i\in \NN $ and $j\in \NN $ with 
$i+j\leq n,\ $ if 
$q_{\sigma ^{i}(x)}^{(j)}(q_{x}^{(i)}(y))=q_{x}^{(i)}(y)$, 
then $|(q_{\sigma ^{i}(x)}^{(j)})'(q_{x}^{(i)}(y))|>1,\ $ 
\end{enumerate}    
we have that 
$\| (q_{x}^{(n)})'(y)\| >M$,\ where $\| \cdot \| $ denotes the 
norm of the derivative with respect to the spherical metric. 
\end{lem}
\begin{proof}
This lemma can be shown by induction on 
$\sharp E$ by means of the same method that was employed in 
Lemma 1.32 in \cite{S4}.
\end{proof}
\begin{lem}
\label{rsporolem0.5}
Let $G=\langle h_{1},\ldots ,h_{m}\rangle $ 
be a finitely generated rational semigroup.
Suppose $\sharp (UH(G)\cap J(G))<\infty $
and $UH(G)\cap J(G)\neq \emptyset .$ 
Then, for each $z\in UH(G)\cap J(G)$, 
there exists an element $g\in G,\ $ 
an element $h\in G$, and a point 
$w\in UH(G)\cap J(G)$ such that 
$h(w)=z,\ $ $g(w)=w$ and $|g'(w)|\leq 1.$  
\end{lem}
\begin{proof}
Suppose that 
there exists a point $z\in UH(G)\cap J(G)$ 
 for which there exists no $(g,h,w)$ in the 
 conclusion of our lemma. 
 Then, by Lemma~\ref{rsporolem0} and 
the Koebe distortion theorem,\ it can 
 easily be seen that, 
 for an arbitrarily small $\epsilon >0$, 
 there exists a positive number $\delta $ 
 and a positive constant $N$ 
 such that, if a point $w_{0}\in UH(G)\cap J(G)$ 
 and an element $g_{0}\in G$ satisfy 
 $g_{0}(w_{0})=z$, then 
 the diameter of the component $V$ 
 of $g_{0}^{-1}(B(z,\delta ))$ containing 
 $w_{0}$  
 is less than $\epsilon $ 
 and $\deg (g_{0}:V\rightarrow B(z,\delta ))\leq N.$ 
 Then, for $\epsilon $ sufficiently small,\ 
 since $G$ is finitely generated 
  and $\sharp (UH(G)\cap J(G))<\infty $, 
  we can easily obtain that 
 there exists a 
 positive constant $N'$ such that,  
 for any element $g_{1}\in G$, 
 and any component $W$ of $g_{1}^{-1}(V),\ $ 
 it follows that $\deg (g_{1}:W\rightarrow V)\leq N'.$ 
 This implies that $z\in SH_{N+N'}(G)$, and 
 this contradicts that 
 $z\in UH(G).$ 
\end{proof}
In the following Lemma~\ref{rsporolem0.6},\ 
Lemma~\ref{rsporolem1},\ and Lemma~\ref{rsporolem2},\ 
we suppose the assumption of Theorem~\ref{rsporosity}.
 Furthermore,\ we assume either (1) $m\geq 2$,\ or 
 (2) $m=1$ and $J(G)\neq \CCI .$  
\begin{lem}
\label{rsporolem0.6}
If $m\geq 2$, then there exists a disk $D$ in $F(G)$ such that 
\begin{enumerate}
\item $\overline{\bigcup _{g\in G}g(D)}\subset F(G)$ and 
\item diam $q_{x}^{(n)}(D)\rightarrow 0$ as 
$n\rightarrow \infty $ uniformly on 
$x\in \Sigma _{m}.$ 
\end{enumerate} 
In particular,\ the fibered rational map $f$ satisfies 
condition (C2).
\end{lem}
\begin{proof}
Suppose the case $m\geq 2.$ 
Let $h\in G$ be an element of degree two or greater.
Since $\emptyset \neq \CCI \setminus \overline{U}\subset F(G)$
 and $UH(G)\cap J(G)\subset U,\ $  
 there exists an attracting periodic point $z_{0}$ of $h$ in 
 $F(G)\setminus U.$ 
 Since $z_{0}\in UH(G)$ and $UH(G)\cap J(G)\subset U$ again,\ 
 it follows that there exists a disk $D$ around $z_{0}$
 such that $\overline{\bigcup _{g\in G}g(D)}\subset F(G).$ 
 Now the second statement of the lemma results from Lemma 1.30 in \cite{S4}.\ 

\end{proof}
\begin{lem}
\label{rsporolem1}
If $UH(G)\cap J(G)\neq \emptyset $, then, 
for each point $z\in UH(G)\cap J(G)$,\   
 there exists an element 
$h\in G$ such that   
$z$ is a parabolic 
fixed point of $h.$ 
Furthermore,\ 
$\{ \beta \in G\mid \beta(z)=z\} 
=\{ g_{z}^{n}\mid n\in \NN \} $ for some 
$g_{z}\in G$ such that $z$ is a parabolic fixed point of $g_{z}.$ 
\end{lem}
\begin{proof}
Let $z\in UH(G)\cap J(G)$ be a point. Then $z\in U.$ 
By Lemma~\ref{rsporolem0.5},\ there exists an element 
$g\in G$,\ an element $\alpha \in G$,\ and 
a point $w\in UH(G)\cap J(G)$ such that 
$\alpha (w)=z,\ g(w)=w$,\ and 
$| g'(w)| \leq 1.$ We write 
$\alpha =h_{i_{k}}\cdots h_{i_{1}}$ and 
$g=h_{j_{l}}\cdots h_{j_{1}}.$ We may assume 
$l\geq k.$ Since $w\in g^{-1}(U)\cap \alpha ^{-1}(U)$,\ 
by the open set condition we obtain 
$(i_{1},\ldots i_{k})=(j_{1},\ldots ,j_{k}).$ Let 
$h:=h_{j_{k}}\cdots h_{j_{1}}\cdot h_{j_{l}}\cdots h_{j_{k+1}}.$ 
Then we easily obtain $h(z)=z.$ Furthermore,\ 
$| h'(z)| =| g'(w)| \leq 1.$   
  
If $\deg (h)=1,\ $ then, by Lemma~\ref{rsporolem0.6},
it follows that $z$ is a repelling fixed point of 
$h.$ However, this is a contradiction. 
If $\deg (h)\geq 2,\ $ then, 
since it is assumed that 
$\sharp (UH(G)\cap J(G))<\infty $,
there exists no recurrent critical point $c$ 
of $h$ in $J(h).$ From \cite{Ma},
 it follows that $z $ is an attracting or 
parabolic fixed point of $h.$ 
 Suppose that $z$ is an attracting fixed point of $h.$ 
 Then there exists an open neighborhood 
 $V$ of $z$ in $U$ such that $h(V)\subset V.$ 
 Let $x\in \Sigma _{m}$  be a point 
 such that $h_{x_{n}} \cdots  h_{x_{1}}=h$, 
 for some $n,\ $ where $x=(x_{1},x_{2},\ldots ).$ 
 Then, by the open set condition,
 it follows that $h_{x'_{n}} \cdots h_{x'_{1}}(V)
 \subset \CCI \setminus U$, 
 for any $(x'_{1},\ldots ,x'_{n})\in \{ 1,\ldots ,m\} ^{n}
 \setminus \{ (x_{1},\ldots ,x_{n})\} .$ 
 Hence, $G$ is normal in $V$, and this is also a 
 contradiction.

 Let $g_{z}$ be an element in 
 $G_{z}:=\{ \beta \in G\mid \beta (z)=z\} $ with 
 minimal word length with respect to the generator system 
 $\{ h_{1},\ldots ,h_{m}\} .$ Let $\beta \in G_{z}.$ 
 Then $z\in \beta ^{-1}(U)\cap g_{z}^{-1}(U).$ 
 By the open set condition,\ it follows that 
 $\beta =\gamma g_{z}$ for some $\gamma \in G.$ 
 Then $\gamma \in G_{z}.$ 
 Applying the same argument again to 
 $\gamma $,\ we see that 
 $\beta =g_{z}^{n}$ for some $n\in \NN .$ 
\end{proof}
\begin{lem}
\label{rsporolem2}
We have that for each 
$(x,y)\in \pi _{\CCI }^{-1}\left(G^{-1}(J(G)\setminus 
UH(G))\right)\cap \tilde{J}(f),\ $  
$\limsup _{n\rightarrow \infty }
d(q_{x}^{(n)}(y),\ UH(G))>0.$
\end{lem}
\begin{proof}
Let $(x,y)$ be a point in 
$\pi _{\CCI }^{-1}\left( G^{-1}\left( J(G)\setminus UH(G)\right) \right) 
\cap \tilde{J}(f).$ 
 Then $q_{x}^{(n)}(y)\in J(G) \setminus UH(G)$, 
for each $n\in \NN .$ 
 
 Assume that $\lim _{n\rightarrow \infty }
 d(q_{x}^{(n)}(y),UH(G))=0.$ We now derive a 
 contradiction.
Let $\epsilon >0$ be a small number and let $A_{\epsilon }$ be the $\epsilon $-neighborhood 
 of $UH(G)\cap J(G)$ in $\CCI .$ 
 Then there exists a number $n_{0}\in \NN $ such that 
 $q_{x}^{(n)}(y)\in A_{\epsilon }$, for each 
 $n\geq n_{0}.$ 
 
  For each $n\geq n_{0},\ $ let $z_{n}\in UH(G)\cap J(G)$ 
  be the unique point such that 
  $d(z_{n},q_{x}^{(n)}(y))<\epsilon .$ 
  Since $g(UH(G))\subset UH(G)$, for each $g\in G$, 
  it may be assumed that 
  $ q_{\sigma ^{n}(x)}^{(1)}(z_{n})=z_{n+1}$, 
  for each $n\geq n_{0}.$ 

  Since $\sharp (UH(G)\cap J(G))<\infty ,\ $ 
  there exists a positive integer $n_{1}\geq n_{0}$ 
  and $l\in \NN $ such that 
  $z_{n_{1}+l}=z_{n_{1}}.$
Let $g_{1}\in G$ be an element such that 
$g_{1}(z_{n_{1}})=z_{n_{1}}.$ 
Let $w\in \{ 1,\ldots ,m\} ^{l}$ be the 
word such that 
$h_{w_{l}}\circ \cdots \circ h_{w_{1}}=g_{1}.$ 
Then, by 
the open set condition, 
$\sigma ^{n_{1}}(x)=w^{\infty }.$ 
Since we assume that $d(q_{x}^{(n)}(y),\ UH(G))\rightarrow 0$ 
as $n\rightarrow \infty ,\ $ 
by $z_{n_{1}+l}=z_{n_{1}}$ it follows that 
$g_{1}^{k}(q_{x}^{(n_{1})}(y))\rightarrow z_{n_{1}}$ 
as $k\rightarrow \infty .$  
Hence, by Lemma~\ref{rsporolem1} it must be that 
$z_{n_{1}}$ is a parabolic fixed point of $g_{1}$, 
and $q_{x}^{(n_{2})}(y)$ belongs to 
$W\cap {\cal P},\ $ where $W$ is a small 
neighborhood of  $z_{n_{1}} $ in $U,\ $ 
${\cal P} $ is the union of  the
attracting petals of $g_{1}$ at $z_{n_{1}}$, 
and  $n_{2}$ is a large positive number with 
$n_{2}\geq n_{1}.$  
Then there exists an open neighborhood $V$ of 
$y$ such that $q_{x}^{(n_{2})}(V)\subset W\cap {\cal P}. $
Taking $W$ sufficiently small and $n_{2}$ sufficiently large,  
it may be assumed that 
$g_{1}^{s}(q_{x}^{(n_{2})}(V))\subset W\cap {\cal P} $, 
for any $s\in \NN .$  
It follows that  
$q_{x}^{(n)}(V)\subset  U$ for each 
$n\in \NN .$ By the open set condition,\ 
for each $n\in \NN $, 
it follows that $q_{x'}^{(n)}(V)\subset \CCI \setminus U$, 
for each $x'\in \Sigma _{m}$ with 
$(x'_{1},\ldots ,x'_{n})\in \{ 1,\ldots ,m\} ^{n}
\setminus \{ (x_{1},\ldots ,x_{n})\} .$
Hence,
$G$ is normal in $V$, and this contradicts that 
$y\in J(G).$      

\end{proof}

We now present a proof of Theorem~\ref{rsporosity}.\\ 
{\bf Proof of Theorem~\ref{rsporosity}.}
By Lemma 2.3 in \cite{S5},\ we have $J(G)\subset \overline{U}.$ 
Suppose $J(G)\neq \overline{U} .$
Then, by Proposition 4.3 in \cite{S4},\ 
 int$J(G)=\emptyset .$
For any $y'\in J(G)$ and $r>0,\ $ we set 
$$h(y',r)=\sup \{ s\mid \exists y''\in \CCI ,\ 
B(y'',s)\subset F(G)\cap B(y',r)\cap U\} $$
and $h(r)=\inf \{ h(y',r)\mid y'\in J(G)\} .$ 
Then, since int $J(G)=\emptyset $ and 
$J(G)\subset \overline{U}$,\  
it follows that $h(r)>0$ for any $r>0.$ 

 Let $\delta _{0}>0$ be a small number. 
 Let $B$ be the $\delta _{0}$-neighborhood of 
 $UH(G)\cap J(G)$ in $\CCI .$     
By Lemma~\ref{keythmlem} and 
Lemma~\ref{rsporolem0.6},\ 
there exists a positive number 
$\delta _{1}$ and a number $N\in \NN $ 
such that, for any $y'\in J(G)\setminus B,\ 
0<\delta \leq \delta _{1}$ and 
any  $V\in c( B(y',2\delta ),\ g),\ $ 
$V$ is simply connected and 
$\deg (g:V\rightarrow B(y' ,2\delta ))\leq N.$ 
By Lemma~\ref{rsporolem1} and the open set condition,\ 
the condition \ref{noshrink} in Theorem~\ref{bundlemeas} 
holds. Hence by 
Lemma~\ref{rsporolem2} and 
Theorem~\ref{bundlemeas},\ we obtain  
\begin{equation}
\label{rsporosityeq1}
\tilde{J}(f)=\bigcup _{x\in\Sigma _{m}}J_{x}.
\end{equation}
Let $y\in J(G)$ be a point. 
Since $\pi _{\CCI }(\tilde{J}(f))=J(G)$
(Proposition 3.2 in \cite{S5}),\  
by (\ref{rsporosityeq1})  
there exists a point $x\in \Sigma _{m}$ 
such that $y\in J_{x}.$ 

 Let $\delta _{2}=\min \{ \delta _{0},\delta _{1}\} .$ 
 Let $r$ be a positive number. We set 
 $B_{n}=q_{x}^{(n)}(B(y,r))$ and 
 $y_{n}=q_{x}^{(n)}(y)$, for each $n\in \NN .$ 
 Since $y\in J_{x},\ $
 there exists a smallest positive integer 
 $n_{0}$ such that 
 diam $B_{n_{0}+1}>\delta _{2}.$ Then there exists a 
 constant $l_{0}$ such that 
 $l_{0}\delta _{2}<$diam $B_{n_{0}}.$ 

 {\bf Case 1.}  $y_{n_{0}}\in J(G)\setminus B.$ 

  By Corollary 2.3 in \cite{Y},\ 
  there exists a constant $K$ that depends only on 
  $N$, and a ball 
  $B(y_{n_{0}},r_{0})\subset B_{n_{0}}$ with 
 $r_{0}\geq $diam $B_{n_{0}}/K\geq \frac{l_{0}\delta _{2}}{K}$ 
such that the element $W\in c( B(y_{n_{0}},r_{0}),\ q_{x}^{(n_{0})})$
containing $y$ is a subset of $B(y,r).$ 
There exists a ball $B(y',\frac{2}{3}h(r_{0}))$ 
included in $B(y_{n_{0}},r_{0})\cap F(G)\cap U.$ 

 Let $D_{0}$ be the element of  
 $c(B(y',\frac{1}{2}h(r_{0})),\ q_{x}^{(n_{0})})$
 contained in $W.$ By the open set condition,\ 
 we obtain that $g^{-1}(U\cap F(G))\subset U\cap F(G)$, for each 
 $g\in G.$ Hence, it follows that $D_{0}\subset F(G)\cap U.$ 
 Let $y''\in D_{0}\cap (q_{x}^{(n_{0})})^{-1}\{ y'\} $ 
 be a point. Then, by Corollaries 1.8 and 1.9 in \cite{S4},\ 
 Dist $(\partial D_{0},y'')\leq M$ for some $M$ 
 that depends only on $N.$ 
 Further, by the same argument that was used in the proof of 
 Theorem~\ref{porositythm}, diam $D_{0}\asymp r.$ 
 Hence, there exists a constant $0<k<1$ that does not 
 depend on $y$ and $r$ such that 
 $B(y'',kr)\subset D_{0}\subset F(G)\cap B(y,r).$ 

 {\bf Case 2.} $y_{n_{0}}\in B.$ 

 Making use of Lemma~\ref{rsporolem1}, the fact  
 that $UH(G)\cap J(G)\subset U,\ $ taking 
 a sufficiently small  $\delta _{0}$,  
 and the method used on pp. 286-287 in \cite{Y}, 
 one can show that there exists a ball $B(y'',k'r)$ 
  in $B(y,r)\cap F(G)$, where $k'$ is a 
  constant with $0<k'<1$ that does not depend 
  on $ y$ and $r.$   
\qed \\  

 We now prove Proposition~\ref{critexp}.\\ 
\noindent {\bf Proof of Proposition~\ref{critexp}.}
 By the open set condition and Lemma 2.3 (f) in \cite{S5},\ 
 $J(G)\subset \overline{U}.$ 
 We now show the following.

 {\bf Claim 1:} There exists an open set $V'$ included in $U\cap F(G)$ 
 such that $h^{-1}(V')\cap V'=\emptyset $, for each $h\in G.$ 

 Before supporting this claim,\ we would like the reader to note that  
the following claim can easily be shown. 

{\bf Claim 2:} If there exists a point $z\in U\cap F(G)$ 
such that $z\in \CCI \setminus \overline{G(z)},\ $ 
then Claim 1 holds with a small open neighborhood 
$V'$ of $z.$ 

 As a proof of Claim 1,\ by the open set condition, 
 we obtain that 
\begin{equation}
\label{critexpeq1}
\bigcup _{j=1}^{m}h_{j}^{-1}(U\cap F(G))\subset U\cap F(G).
\end{equation} 
Suppose the equality does not hold in (\ref{critexpeq1}).
Then there exists a point $z\in U\cap F(G)$ such that 
$h_{j}(z)\in \CCI \setminus U$ for each $j=1,\ldots ,m.$ 
Hence, by the open set condition,\  
$z\in \CCI \setminus \overline{G(z)}.$ By Claim 2,\ 
Claim 1 holds.

 Hence, it may be assumed that 
 \begin{equation}
 \label{critexpeq2}
\bigcup _{j=1}^{m}h_{j}^{-1}(U\cap F(G))= U\cap F(G).
  \end{equation}
 Let $\alpha :U\cap F(G)\rightarrow U\cap F(G)$ be the 
 map defined as: $\alpha (z)=h_{j}(z)$ if $z\in h_{j}^{-1}(U\cap F(G)).$ 
 This is well defined by (\ref{critexpeq2}) and the open set 
 condition.

 Let $z\in U\cap F(G)$ be a point. If $z\in \CCI \setminus \overline{G(z)},\ $
 then, by Claim 2, Claim 1 follows. Hence, it may be assumed that 
 $z\in \overline{G(z)}$; i.e., 
 \begin{equation}
 \label{critexpeq3}
 z\in \overline{\bigcup _{n=0}^{\infty }\{ \alpha ^{n}(z)\}}.
 \end{equation}
Let $W$ be the connected component of $U\cap F(G)$ containing 
$z.$ By (\ref{critexpeq3}), there exists a 
smallest positive integer $n$ with $\alpha ^{n}(W)\subset W.$ 
By (\ref{critexpeq3}) and the open set condition,\ 
one of the following two
cases holds.

 Case 1: $W$ is included in an attracting or a parabolic basin of an element
  $g\in G,\ $ 
  $z$ is a fixed point in the basin, and 
  $g|_{W}=\alpha ^{n}|_{W}.$ 

 Case 2: $W$ is included in a Siegel disk or a Herman ring 
 of an element $g\in G$ of degree $2$ or greater, and 
 $g|_{W}=\alpha ^{n}|_{W}.$

 If Case 1 holds,\ then there exists an open set 
 $V'$ included in $W$ with $\alpha ^{-l}(V')\cap V'=\emptyset $, 
 for each $l\in \NN $; i.e., $h^{-1}(V')\cap V'=\emptyset $, 
 for each $h\in G.$ 

 If Case 2 holds,\ then, for $V'$ in a connected component 
 $A$ of $\alpha ^{-n}(W)$ with $A\cap W=\emptyset ,\ $ it follows that 
 $\alpha ^{-l}(V')\cap V'=\emptyset $, for each $l\in \NN $; i.e., 
 $h^{-1}(V')\cap V'=\emptyset $, for each $h\in G.$ 

 Hence, we have proved Claim 1. Let $V'$ be an 
 open set included in $U\cap F(G)$ such that 
 $h^{-1}(V')\cap V'=\emptyset $ for each $h\in G.$ 
 Then, by the open set condition, 
 $g^{-1}(V')\cap h^{-1}(V')=\emptyset ,\ $ if $
 g,h\in G$ and $g\neq h.$ Furthermore, the postcritical set 
  $ P(G)$
 of $G$,\ which equals 
$\overline{\bigcup _{g\in G}g\left(\bigcup _{j=1}^{m}
\{ \mbox{critical values of }h_{j}\} \right) }$,\  
  does not accumulate in $V'.$ Let $V$ be an open disk included in 
 $V'\setminus P(G).$ Then, it follows that 
 $ \int _{V}\limits\sum _{h\in G}\limits\sum _{\alpha }
 \| \alpha '(z)\| ^{2}\ dm(z) <\infty ,\ $
 where $\alpha $ runs over all the well-defined inverse branches 
 of $h$ on $V.$ Hence, for almost every $x\in V$ with respect to 
 the Lebesgue measure,\ it follows that $S(2,x)<\infty .$   
 \qed \\ 

 We now prove Theorem~\ref{uhfindim}.
 We first present several lemmas.
\begin{lem}
\label{forhdimlem2}
Let $G$ be a rational semigroup. Assume that 
$\infty \in F(G) $ and that, for each $x\in E(G)$, there exists an element 
$g\in G$ such that $g(x)=x$ and $|g'(x)|< 1.$ 
Let $A$ be a subset of $J(G).$ 
Suppose that there exist positive constants $a_{1}$
$,a_{2}$ and $c$ with $0<c<1$ such that, for each 
$x\in A,\ $ 
there exist two sequences $(r_{n})$ and $(R_{n})$ of positive real numbers 
and a sequence $(g_{n})$ of the elements of $G$ satisfying all of the following conditions:
\begin{enumerate}
\item \label{forhdim21} $r_{n}\rightarrow 0 $  and, for each 
$n,\ $ $0<\frac{r_{n}}{R_{n}}<c$ and $g_{n}(x)\in J(G).$     
\item \label{forhdim22} for each $n,\  g_{n}(D(x,R_{n}))\subset D(g_{n}(x),a_{1}).$
\item \label{forhdim23}for each $n\ g_{n}(D(x,r_{n}))\supset D(g_{n}(x),a_{2}).$ 
\end{enumerate}
Then, 
$ \mbox{dim}_{H}(A)\leq s(G).$ 
\end{lem}
\begin{proof} 
We may assume that $\sharp (J(G))\geq 3.$ 
Let $\delta \geq s(G)$ be a number and $\mu $ a 
$\delta $-subconformal measure. 
Using the method employed in the proof of Lemma 5.5 in \cite{S4},\ 
one can show that there exists a constant $c'>0$ 
that does not depend on $n\in \NN $ and $x\in A$ such that 
$ \mu (D(x,r_{n}))\geq c'r_{n}^{\delta }.$ 
From this and Theorem 7.2 in \cite{Pe},\ 
it follows that dim$_{H}A\leq \delta .$   
\end{proof}
\begin{df}({\bf Conical set for rational semigroups})
Let $G$ be a rational semigroup. Let $N\in \NN $ and 
$r>0.$ We denote by $J_{con}(G,N,r)$ the set of points 
$z\in J(G)$ such that, for any $\epsilon >0,\ $ 
there exists an element $g\in G$ such that $g(z)\in J(G)$ and 
the element $U\in c(B(g(z),r),\ g)$ containing $z$ satisfies 
the following conditions:
\begin{enumerate}
\item diam $U\leq \epsilon $,\ 
\item $U$ is simply connected,\ and 
\item deg$(g:U\rightarrow B(g(z),r))\leq N.$  

\end{enumerate}
We set $J_{con}(G,N)=\cup _{r>0}J_{con}(G,N,r)$ and 
$J_{con}(G)=\cup _{N\in \NN }J_{con}(G,N).$ 
\end{df}
\begin{prop}
\label{condim}
Let $G$ be a rational semigroup. Assume that 
$F(G)\neq \emptyset $ and 
that, for each $x\in E(G),\ $ there exists an element 
$g\in G$ such that $g(x)=x$ and $|g'(x)|<1.$ 
Then, 
$\mbox{dim}_{H}(J_{con}(G))\leq s(G).$  
\end{prop}
\begin{proof}
We now have only to show the following: 

 {\bf Claim:} For fixed  $N\in \NN $ and $r>0,\ $ 
 dim$_{H}(J_{con}(G,N,r))\leq s(G).$

 We may assume that $\infty \in F(G).$ 
 Let $x\in J_{con}(G,N,r)$ be a point. 
 Then there exists a sequence $(g_{n})$ in $G$ 
 such that, for each $n\in \NN $, we have that $g_{n}(x)\in J(G),\ $ 
 $ \deg (g_{n}:V_{n}(r)\rightarrow D(g_{n}(x),r))\leq N , $  
 $V_{n}(r)$ is simply connected and 
 diam $V_{n}(r)\rightarrow 0$ as $n\rightarrow \infty ,\ $ 
 where $V_{n}(r)$ is the element of $c(D(g_{n}(x),r),\ g_{n})$ 
 containing $x.$ 
  Let $\varphi _{n}:D(0,1)\rightarrow V_{n}(r)$ be the 
  Riemann map such that $\varphi _{n}(0)=x.$ By the 
  Koebe distortion theorem, for each $n,\ $ 
  $V_{n}(r)\supset D(x,\frac{1}{4}|\varphi _{n}'(0)|).$ 
  By Lemma~\ref{ylem1} and the Koebe distortion theorem,\ 
  there exists an $\epsilon >0$ such that, 
  for each $n\in \NN ,\ $ 
  $V_{n}(\epsilon r)\subset D(x,\frac{1}{8}|\varphi _{n}'(0)|).$ 
  Since diam $V_{n}(r)\rightarrow 0$ as $n\rightarrow \infty ,\ $ 
  it follows that $|\varphi _{n}'(0)|\rightarrow 0$ as $n\rightarrow \infty .$ 
  The application of Lemma~\ref{forhdimlem2}\ proves the claim.       
\end{proof}
\begin{df}({\bf Good points for finitely generated rational semigroups})
\label{Good pts in J}
Let $G=\langle h_{1},\ldots ,h_{m}\rangle $ 
be a rational semigroup. 
Let $f: \Sigma _{m}\times \CCI \rightarrow 
\Sigma _{m}\times \CCI $ be the fibered rational 
map associated with the generator system 
$\{ h_{1},\ldots ,h_{m}\} .$ 
Then we set 
$J_{good }(G)=\pi _{\CCI }(\tilde{J}_{good}(f)).$ 
Note that this definition 
does not depend on the choice of 
any generator system of $G$ that consists of 
finitely many elements.
\end{df}
We now prove the following theorem.
\begin{thm}
\label{jgooddim}
Let $G=\langle h_{1},\ldots ,h_{m}\rangle $ be a finitely 
generated rational semigroup with $F(G)\neq \emptyset.$ 
Let $f:Y\rightarrow Y$ be the fibered rational map 
associated with the generator system 
$\{ h_{1},\ldots ,h_{m}\} ,\ $ where $Y=\Sigma _{m}\times \CCI .$ 
Suppose that $f$ satisfies condition (C1) and 
that, for each $x\in \Sigma _{m},\ $ 
the boundary of $\hat{J}_{x}(f)\cap UH(f)$ in 
$Y_{x}$ does not separate points in $Y_{x}.$ Then, 
$J_{good}(G)\subset J_{con}(G)$ and 
$ \dim _{H}(J_{good}(G))\leq s(G)\leq s_{0}(G).$   
\end{thm}
\begin{proof}
We may assume that $\sharp (J(G))\geq 3.$ 
We now show the following:

 {\bf Claim:} If $E(G)\neq \emptyset ,\ $ 
 then, for each $x\in E(G)$, there exists an 
 element $g\in G$ such that $g(x)=x$ and $|g'(x)|<1.$ 

 If there exists an element $h\in G$ with $\deg (h)\geq 2,\ $ 
 then this claim is trivial. Suppose that each element of $G$ 
 is of degree one. By Lemma 2.3 (d) in \cite{S5},\  
 $\sharp (E(G))\leq 2.$ 
Since $f$ satisfies condition (C1),\ 
for each $i,\ h_{i}$ is loxodromic. 
Since $h_{i}(E(G))=E(G)$, for each $i,\ $ 
it must be that each $x\in E(G)$ is fixed by 
$h_{i}$, for each $i.$  Let $x\in E(G)$ be a point.
Suppose $|h_{i}'(x)|>1$, for each $i.$ 
 Then, $J(G)=\{ x\} $, and this is a contradiction, 
 since we assume that $\sharp (J(G))\geq 3.$ 
 Hence, $|h_{i}'(x)|<1$ for some $i.$ Hence, the claim holds.

 The statement of this theorem follows from the claim,\ 
 the second statement in 
Proposition~\ref{measprop2},\ 
Proposition~\ref{condim} and 
Theorem 4.2 in \cite{S2}.
\end{proof}
We now prove Theorem~\ref{uhfindim}.\\ 
\noindent {\bf Proof of Theorem~\ref{uhfindim}.}
This follows from Lemma~\ref{rsporolem0.6},\ 
$\pi _{\CCI }(\tilde{J}(f))=J(G)$ (Proposition 3.2 in \cite{S5}),\ 
 Lemma~\ref{rsporolem2} 
and Theorem~\ref{jgooddim}, when $m\geq 2.$ 

 We now suppose $m=1$ and $J(G)\neq \CCI .$ 
 Then, by \cite{Ma}, $J_{good}(G)\subset J_{con}(G).$ 
 By Lemma~\ref{rsporolem2} and Proposition~\ref{condim}, 
 it follows that $\dim _{H}(J(G))\leq s(G).$ 
Moreover, by Theorem 4.2 in \cite{S2}, we have 
$s(G)\leq s_{0}(G).$ 

 Hence, we have proved Theorem~\ref{uhfindim}.   
\qed

\end{document}